\documentclass[12pt,titlepage]{article}

\renewenvironment{abstract}{\null\vfil
\begin{center}
  \bfseries \abstractname
\end{center}\begin{quote}}{\end{quote}\vfil\null}

\usepackage[tbtags]{amsmath}
\usepackage{amsopn,amsthm,amsfonts,amssymb}

\numberwithin{equation}{section}

\DeclareMathOperator{\Res}{Res}

\makeatletter
\renewcommand{\subsubsection}{\@startsection
{subsubsection}
{3}
{0mm}
{\baselineskip}
{-0.5\baselineskip}
{\normalfont\normalsize\bfseries}}
\makeatother

\newtheorem{theorem}{Theorem}[section]
\newtheorem{lemma}[theorem]{Lemma}
\newtheorem{proposition}[theorem]{Proposition}

\newtheorem{corollary}[theorem]{Corollary}

\theoremstyle{definition}

\theoremstyle{remark}

\newtheorem*{acknow}{Acknowledgments}

\begin{document}

 \title{Exact operator solution of the Calogero-Sutherland model}
\author{Luc~\textsc{Lapointe} and Luc~\textsc{Vinet}\\
\small
\begin{tabular}{c}
Centre de recherches math{\'e}matiques\\
Universit{\'e} de Montr{\'e}al\\
C. P.~6128,  succ.~Centre-ville\\
Montr{\'e}al QC H3C 3J7
\end{tabular}}
\date{\textbf{CRM-2272}\\[\bigskipamount]
  April 1995}
\maketitle

\begin{titlepage}

\begin{abstract}
The wave functions of the Calogero-Sutherland model are known to be
expressible in terms of Jack polynomials.  A formula which allows to obtain
the wave functions of the excited states by acting with a string of creation
operators on the wave function of the ground state is presented and derived.
The creation operators that enter in this formula of Rodrigues-type for the
Jack polynomials  involve Dunkl operators.
\end{abstract}

\bigskip

\textbf{Short title:} Solution of the CS Model

\bigskip


\end{titlepage}

\section{Introduction}
Exactly solvable models are of great help in the understanding of quantum
many-body physics.  The Calogero-Sutherland (CS) \cite{1,2,3} model, which
describes a
system of $N$ particles on a circle interacting pairwise through long range
potentials, is generating a lot of attention in this connection, in particular
because it provides a fully solvable model in which the ideas of fractional
statistics can be tested \cite{4}.  There is thus considerable interest in
identifying
the algebraic structure responsible for the solvability of this model.

The spectrum of the CS Hamiltonian can be interpreted as the energy of a
collection of free quasi-particles obeying a generalized exclusion principle.
Recent computations \cite{5,6,7} of some correlation functions have confirmed
this point
of view and shown that the exclusion statistics of quasi-particles and
quasi-holes is consistent with the anyon statistics of the real particles.
The calculation of these quantities proved possible because of the following
circumstance:  the eigenfunctions of the CS Hamiltonian are given in terms of
Jack polynomials \cite{8,9,10}.  These polynomials form a basis for the ring of
symmetric
functions and enjoy algebraic properties that allow to carry out analytically
the computation of various dynamical functions of the CS model.  These
polynomials also appear in related areas like the characterization of the
Virasoro and $W_N$ algebras singular vectors \cite{11,12} and in the
construction of Yangian
modules \cite{13,14}.

 We will show in this paper that the wave functions of the CS Hamiltonian and
hence the Jack polynomials can be obtained by applying a string of creation
operators on the ground state wave function.  After reviewing in section~2,
the spectrum and the eigenstates of the CS Hamiltonian, we shall present and
discuss this operator solution of the CS model in section~3.  Proofs will be
deferred mainly to section~4.  Conclusion and outlook will form the content of
section 5.

\section{Spectrum and eigenstates of the  Hamiltonian}
The CS model describes a system of $N$ particles on a circle.  We shall denote
by $L$ the perimeter of that circle and by $x_i$, $i=1, \dots , N$; $0 \leq x_i
\leq L $, the positions of the particles.

\subsection{Hamiltonian and ground state}
The quantum Hamiltonian $H_{CS}$ is
\begin{equation}
H_{CS}
	= - \sum_{j=1}^N \frac{\partial^2}{\partial x_j^2} +2 \beta (\beta-1)
\sum_{j<k} \frac{1}{d^2(x_j-x_k)},
\end{equation}
where $\beta$ is a constant and
\begin{equation}
d(x_j-x_k)
	= \frac{L}{\pi} \sin \frac{\pi}{L} (x_j-x_k),
\end{equation}
is the chord length between the positions of the particles $j$ and $k$.  Note
the symmetry under the exchange of $\beta$ and $1-\beta$.  For $\beta$ real,
that is $\beta (\beta-1) \geq -1/4$, $H_{CS}$ is known to be stable and to have
no bound states.  The momentum operator $P_{CS}$ is $P = -i \sum_{j=1}^N
{\partial}/{\partial x_j}$.  The operators $P$ and $H_{CS}$ are self-adjoint
with respect to the inner product $(f,g)= \int_0^L  dx_1 \dots \int_0^L \,dx_N
f(x)
\overline {g(x)}$.  Moreover, since $H_{CS}$ can be written in the form
\cite{12}
\begin{equation}
H_{CS}
	= \sum_{j=1}^N A_j^{+}(\beta) A_j(\beta) + E_{0 },
\end{equation}
where
\begin{equation}
A_j(\beta)
	= -i \frac{\partial}{\partial x_j} + i \frac{\pi}{L} \beta \sum_{k \neq j}
\cot \bigl[\frac{\pi}{L} (x_j-x_k) \bigr],
\end{equation}
we see that $H_{CS}$ is bounded from below with
\begin{equation}
E_{0}
	= \frac{1}{3} \Bigl(\frac{\pi}{L}\Bigr)^2 \beta^2 N(N^2-1),
\end{equation}
the ground state energy.  There are a priopri two wavefunctions satisfying
$H_{CS} \psi_{0} = E_{0} \psi_{0}$, one being annihilated by the operators
$A_j(\beta)$ and the other by the operators $A_j (1-\beta)$.  Only the first
is normalizable however for all values of $\beta$.  It is of Jastrow type and
given explicitely by
\begin{equation}
\psi_{0} (x)
	= \prod_{j<k} \sin \Bigl[ \frac{\pi}{L} (x_j-x_k) \Bigr] ^{\beta}.
\end{equation}
The coupling $\beta$ controls the statistics of the real particles in the
ground state.

\subsection{Excited states}
The wave functions of the excited states are written in the form $\psi (x)=
\phi (x) \psi_{0} (x)$ where $\phi(x)$ is required to be symmetric in order
for $\psi$ to behave like $\psi_{0}$ under the exchange of particles.  It
proves
convenient to use the variables
\begin{equation}
z_j
	= e^{2 \pi i x_j/L}.
\end{equation}
In these coordinates, the ground state wave function is (up to a constant)
$\Delta^{\beta}(z) = \prod_{j<k} (z_j-z_k)^\beta \prod_{k} z_k^{-\beta
(N-1)/2}$
and the Schrodinger equation $H_{CS} \psi = E \psi$ is transformed into the
following equation for $\phi$ :
\begin{equation}
H \phi
	= \Bigl(\frac{L}{2 \pi}\Bigr)^2 (E-E_{0}) \phi,
\end{equation}
where
\begin{equation}
H
	= \Bigl(\frac{L}{2 \pi}\Bigr)^2 \Delta^{-\beta} H_{CS}~\Delta^{\beta} =
H_1+\beta H_2,
\end{equation}
\begin{subequations}
\begin{align}
&H_1
	=\sum_{j=1}^N \Bigl( z_j \frac{\partial}{\partial z_j} \Bigr)^2,
\\
&H_2
	= \sum_{j<k} \Bigl( \frac{z_j+z_k}{z_j-z_k} \Bigr) \Bigl(
z_j \frac{\partial}{\partial z_j}-z_k \frac{\partial}{\partial z_k} \Bigr).
\end{align}
\end{subequations}
The momentum operator becomes $P= \Delta^{-\beta} P \Delta^{\beta}={2 \pi}/{L}
\sum_{j} z_j  {\partial}/{\partial z_j}$ and commutes with $H$.  We thus
supplement (2.8) with
\begin{equation}
P \phi
	= \kappa~\phi .
\end{equation}
Let $\phi^{\prime}=G^{q}\phi$ with
\begin{equation}
G
	= \prod_{j=1}^N z_j, \qquad q \in R.
\end{equation}
If $\phi$ obeys (2.8) and (2.11), it is easy to see that $\phi^{\prime}$ will
also be an eigenfunction of $H$ and $P$ with eigenvalues $ ({L}/{2 \pi})^2
(E-E_{0})+2Nq({L}/{2 \pi}) \kappa +(Nq)^2$ and $\kappa+Nq$ respectively.
Multiplication
by $G$ thus implements Galilei boosts.  Finally, one readily notices that $H$
is $S_N$-invariant.  This implies that the space of symmetric functions of
degree $n \in \mathbb N$ is stable under the action of $H$.  In fact, it is
shown
\cite{5,6,7} that
apart from factors of the form (2.12), the eigenfunctions $\phi$ of $H$ and $P$
are given in terms of a specific basis for these symmetric functions which is
known as that of the Jack polynomials $J_{\lambda}(z;1/\beta)$.

\subsection{Jack polynomials}
These  polynomials are labelled by partitions $\lambda$ of their degree $n$,
that
is sequences $\lambda = (\lambda_1,\lambda_2, \dots)$ of non-negative integers
in
decreasing order $\lambda_1 \geq \lambda_2 \geq \dots$ such that $n=\lambda_1+
\lambda_2 + \dots$   Let $\lambda$ and $\mu$ be two partitions of $n$.  In the
dominance ordering, we have $ \lambda \geq \mu$ if $\lambda_1+\lambda_2+\dots
+
\lambda_i \geq \mu_1 + \mu_2 +\dots +\mu_i$, for all $i$.  Two natural bases
for
the space of symmetric functions are conventionally used to define the Jack
polynomials:
\begin{enumerate}
\item[(i)] the power sum symmetric functions $p_{\lambda}$ which in terms of
the
 power sums
\begin{equation}
p_i
	=\sum_k z_k^i,
\end{equation}
are given by
\begin{equation}
p_{\lambda}
	=p_{\lambda_1}p_{\lambda_2}\dots ,
\end{equation}
\item[(ii)]  the monomial symmetric functions $m_{\lambda}$ which are
\begin{equation}
m_{\lambda}
	= \sum_{\text{distinct permutations}}z_1^{\lambda_1} z_2^{\lambda_2}\dots
\end{equation}
\end{enumerate}
To the partition $\lambda$ with $m_i$ parts equal to $i$, we associate the
number
\begin{equation}
z_{\lambda}
	= 1^{m_1} m_1!~ 2^{m_2} m_2! \dots
\end{equation}
We then introduce the following scalar product on the space of symmetric
functions
\begin{equation}
\langle p_{\lambda}, p_{\mu} \rangle
	= \delta_{\lambda,\mu} z_{\lambda}\beta^{-l(\lambda)},
\end{equation}
where $l(\lambda)$ is the number of parts of $\lambda$.  The Jack polynomials
$J_{\lambda}(z_1,\dots,z_N;1/\beta)$ are then uniquely defined \cite{8,9,10} as
the symmetric polynomials
satisfying the following two conditions:
\begin{subequations}
\begin{align}
\langle & J_{\lambda},J_{\mu}\rangle
	= 0  \quad \text{ if } \lambda \neq \mu , \\
J_{\lambda}\Bigl(z; & 1/\beta \Bigr)
	= m_{\lambda}+\sum_{\mu < \lambda} v_{ \lambda \mu}(\beta)
m_{\mu}.
\end{align}
\end{subequations}
While no explicit formula has yet been obtained for these polynomials, they
have
been shown to obey a number of interesting properties.  For instance,  they are
also orthogonal under the norm \cite{10}
\begin{equation}
(\phi_1,\phi_2)
	= \int \frac{d\theta_1}{2 \pi} \dots \frac{d\theta_N}{2 \pi}
\prod_{j<k} \bigl|(z_j-z_k)\bigr|^{2 \beta} \phi_1(z)\overline {\phi_2(z)},
\qquad z_j=
e^{i \theta_j},
\end{equation}
which is induced from the scalar product associated to the original quantum
mechanical problem $(\theta_i=2 \pi x_i/L)$.  We can therefore replace $
\langle J_{\lambda},J_{\mu}\rangle$ by $(J_{\lambda},J_{\mu})$ in (2.18a) to
define
the Jack functions, these two scalar product being proportional.  Moreover,
when $l(\lambda) \leq N$, the Jack polynomials $J_{\lambda}(z_1, \dots, z_N;
1/\beta)$ are shown \cite{8,9,10} to obey a differential equation which
coincides
with (2.8)--(2.10).

\subsection{Diagonalization of $H$}

This last result can be presented following Sutherland's argument \cite{3} on
how to
triangulate the CS Hamiltonian.  Since the ground state incorporates
correlations, it is reasonable to expect, when writing the wave functions in
the
form $\psi=\phi ~\psi_{0}$, that $\phi$ will look much like the wave functions
of free particles which are symmetric monomials
$m_{\lambda}=\sum_{\text{perm}}\prod_j
z_j^{\lambda_j}$. It is easy to see that $H=H_1+\beta H_2$ is self-adjoint with
respect to the scalar product (2.19).  It thus has a set of orthogonal
eigenfunctions which form a basis for the symmetric functions.  Consider now
the action of $H$ on $m_{\lambda}$.  One readily finds, using (2.10) that
\begin{equation}
\begin{split}
&H_1 m_{\lambda}
	= \biggl(\sum_j \lambda_j^2\biggr)  m_{\lambda},
\\
&H_2 m_{\lambda}= \varepsilon_2 m_{\lambda}+ \sum_{\mu < \lambda} c_{\lambda
\mu}
m_{\mu},
\end{split}
\end{equation}
where $c_{\lambda \mu}$ are some coefficients and
\begin{equation}
\varepsilon_2
	= \sum_{j<k}(\lambda_j -\lambda_k)=\sum_{j}(N+1-2j)\lambda_j.
\end{equation}
We thus find that $H$ is triangular in the symmetric monomial basis.  Its
eigenvalues are
thus
\begin{equation}
\Bigl(\frac{L}{2 \pi}\Bigr)^2 (E-E_{0})
	= \sum_{j=1}^N \bigl(\lambda_j^2+\beta(N+1-2j)
\lambda_j \bigr).
\end{equation}
The corresponding eigenfunctions $\phi_{\lambda}$ will have an expansion of the
form $\phi_{\lambda}=  m_{\lambda}+ \sum_{\mu < \lambda} v_{\lambda \mu}
m_{\mu}$ and satisfy $(\phi_{\lambda},\phi_{\mu})=0$.  Since these two
conditions uniquely define the Jack polynomials, we must have $\phi_{\lambda}
(z)=J_{\lambda}(z;1/\beta)$.

To sum up, allowing for the Galilean boosts, the eigenfunctions of the
Calogero-Sutherland Hamiltonian are
\begin{equation}
\psi_{\lambda,q}(z)
	=\biggl( \prod_{i=1}^N z_i \biggr)^{q-(N-1)\beta/2}
\prod_{i<j}(z_i-z_j)^{\beta}~J_{\lambda}(z;1/\beta),
\end{equation}
where $q$ is an arbitrary real number and $\lambda$ ranges over all partitions
of all non-negative integers such that $l(\lambda) \leq N-1$.  The number of
parts in $\lambda$ is restricted to be strictly smaller than the number of
particles in order to avoid double counting.  Indeed, it is immediate to
convince oneself that the solution $\bigl( \prod_{i=1}^N z_i \bigr)
J_{\lambda}(z;1/\beta)$ of (2.18) is actually the solution $J_{\lambda+1}
(z;1/\beta)$ where $\lambda+1=(\lambda_1+1,\lambda_2+1, \dots)$  We therefore
set
$\lambda_N=0$ absorbing any value that this part might take in the arbitrary
Galilei transformation parametrized by $q$.

The eigenvalues $\kappa_{\lambda}$ and $E_{\lambda}$ that $P$ and $H_{CS}$ have
when acting on the wavefunctions (2.23) can be nicely presented if one
introduces the quantities
\begin{equation}
\kappa_i
	= \frac{2 \pi}{L} \bigl[ \lambda_i + \beta(N+1-2i)+q\bigr].
\end{equation}
In terms of these
\begin{gather}
\kappa_{\lambda}
	=\frac{2 \pi}{L} (n+Nq)=\sum_{i=1}^N \kappa_i, \\
E_{\lambda}
	=\sum_{i=1}^N \kappa_i^2.
\end{gather}
The spectrum of the CS model is that of free quasi-particles with
quasi-momenta $\kappa_i$.  The neighboring quasi-momenta satisfy
$\kappa_i-\kappa_{i+1} \geq 2 \pi \beta/L$
which indicates that the quantum excitations of the CS system
obey a generalized exclusion statistics.

\section{Operator solution}
  We shall now present a formula that gives the eigenfunctions of the CS
Hamiltonian and hence the Jack polynomials, through the action of a string of
creation operators on the ground state wave function.  The proof of this
formula will be given in the next section.

\subsection{Dunkl operators}
The creation operators  that will provide this operator
solution of the CS model will be constructed in terms of the so-called Dunkl
 operators  $\nabla_i$ \cite{15}.  These operators are defined as follows
\begin{equation}
\nabla_i
	= \frac{\partial}{\partial z_i} + \beta \sum_{\begin{subarray}{c} j=1\\ j \ne
i\end{subarray}}^N \frac{1}
{(z_i-z_j)}(1-K_{ij}),
\end{equation}
where $K_{ij}=K_{ji}$, $K_{ij}^2=1$, is the operator that permutes the
variables $z_i$ and $z_j$:
\begin{equation}
K_{ij} z_j
	=z_i  K_{ij} .
\end{equation}
The Dunkl operators are easily found to have the following properties:
\begin{subequations}
\begin{align}
&~ [\nabla_i,\nabla_j]
	=0,\\
& K_{ij}\nabla_j
	=\nabla_i K_{ij},\\
[\nabla_i,z_j]=& \delta_{ij}
\biggl(1+\beta \sum_{l=1}^N K_{il}\biggr) - \beta K_{ij}.
\end{align}
\end{subequations}
We now define
\begin{equation}
D_i
	\equiv z_i \nabla_i .
\end{equation}
In terms of these operators, the operator $H$ given in (2.9) and (2.10) takes a
remarkably simple form:
\begin{equation}
H
	= \Res  \sum_{i=1}^N D_i^2,
\end{equation}
where $\Res   X$ means that the action of $X$ is restricted to symmetric
functions
of the variables $z_1, \dots ,z_N$.  Actually, $H$ belongs to a set of mutually
commuting operators.  Being completely integrable, the CS system admits
$N$ functionally independant constants of motion that are in involution
\cite{16}.  Modulo
conjugation by $\Delta^{\beta}$ (see (2.9)), these are
\begin{equation}
\begin{split}
L_j
	= & \Res  \sum_{i=1}^N (D_i)^j, \qquad j=1, \dots ,N, \\
& [L_k,L_j]=0, \qquad L_2=H.
\end{split}
\end{equation}

\subsection{Main result}
We need additional notation to present our main result.  Let
$J =\{ j_1,j_2,\dots , j_{\ell} \}$ be  a set of cardinality $|J|=\ell$ made of
integers
such that
\begin{equation}
j_{\kappa} \in \{1,2,\dots,N\}, \qquad  j_1<j_2<\dots<j_{\ell}.
\end{equation}
Introduce the operators $D_{k, J}$ labelled by these sets $J$ and by
non-negative integers $k$
\begin{equation}
D_{k,J}
	=\Bigl(D_{j_1}+k\beta\Bigr)\Bigl(D_{j_2}+(k+1)\beta\Bigr)
\dots \Bigl(D_{j_{\ell}}+(k+\ell-1)\beta\Bigr).
\end{equation}
We now set
\begin{equation}
B_{i,J}^{+} = \sum_{J^{\prime} \subset J;|J^{\prime}|=i} z_{J^{\prime}}
D_{1,J^{\prime}}, \qquad i<N ,
\end{equation}
with
\begin{equation}
z_{J^{\prime}}
	=\prod_{i \in J^{\prime}} z_i.
\end{equation}
The sum in (3.9) is over all subsets $J^{\prime}$ of $J$ that are of
cardinality
$i$.  When $J=\{1,\dots,N \}$, we shall write for short
\begin{equation}
B_i^{+}
	\equiv B_{i, \{1,\dots,N\}}^{+}, \qquad i<N.
\end{equation}
We shall also take
\begin{equation}
B_N^{+}
	= z_1 z_2 \dots z_N.
\end{equation}
This last operator is identical to the Galilei boost operator of (2.12).

We can now state the following result.

\begin{theorem}  The Jack polynomials $J_{\lambda}(z;1/\beta)$ associated to
partitions $\lambda=(\lambda_1,\lambda_2, \dots ,\lambda_{N-1})$ are given by
\begin{equation}
J_{\lambda}(z;1/\beta)
	=c_{\lambda}^{-1} \bigl( B_{N-1}^{+}\bigr)^{\lambda_{N-1}} \dots
\bigl( B_{2}^{+}\bigr)^{\lambda_2-\lambda_3}\bigl( B_{1}^{+}\bigr)^
{\lambda_{1}-\lambda_2} \cdot 1,
\end{equation}
with the constant $c_{\lambda}$ equal to
\begin{equation}
c_{\lambda}
	= \prod_{k=1}^{N-1} c_k(\lambda_1,\dots,\lambda_{k+1};\beta),
\end{equation}
where $\lambda_N \equiv 0$ and
\begin{multline}
c_k(\lambda_1,\dots,\lambda_{k+1};\beta)
	= (\beta)_{\lambda_k-\lambda_{k+1}} (2 \beta+\lambda_{k-1}-
\lambda_{k})_{\lambda_k-\lambda_{k+1}} \dots\\
\dots (k \beta + \lambda_1-\lambda_{k})_{\lambda_k-\lambda_{k+1}}.
\end{multline}
\end{theorem}
\noindent In (3.15), $(\beta)_n$ stands for the Pochhammer symbol that is
$(\beta)_n$ = $\beta(\beta+1) \dots (\beta+n-1)$, $(\beta)_0 \equiv 1$.

We can of course conjugate the
creation operators with the ground state wave function $\psi_{0}(z)=
\bigl( \prod_{i=1}^N z_i \bigr) ^{-(N-1)\beta/2} \prod_{i<j} (z_i-z_j)^{\beta}$
and write
\begin{equation}
{\widetilde B}_i^{+} = \psi_{0} B_i^{+} \psi_{0}^{-1}.
\end{equation}
In view of (2.23) and (3.13), we then have the following formula.

\begin{corollary}  The eigenfunctions of the CS Hamiltonian are given by
\begin{equation}
\psi_{\lambda,q}(z)
	= c_{\lambda}^{-1} ({\widetilde B}_N^{+})^q
({\widetilde B}_{N-1}^{+})^{\lambda_{N-1}} \dots
({\widetilde B}_{2}^{+})^{\lambda_2-\lambda_3} ({\widetilde B}_{1}^{+}
)^{\lambda_{1}-\lambda_2} \psi_{0}(z) .
\end{equation}
\end{corollary}
This shows how the wave functions $\psi_{\lambda,q}$ of the excited states can
be obtained by applying iteratively creation operators on the ground state
wave function of the system.

\subsection{Remarks}
Formulas (3.13) and (3.14) will be proved in the next section.  Let us here
comment on some of their features.

\subsubsection{} It is not immediately obvious that the right-hand side of
(3.13) yields
symmetric polynomials of degree $n=\lambda_1+ \dots +\lambda_{N-1}$ in the
variables
$z_1,\dots ,z_N$ . This however is easily seen to follow from the properties of
the Dunkl operators.  We shall often use the notation $\Res  ^{\{i,j,k, \dots
\}}X$
or $\Res  ^{J} X$ to indicate that $X$ is taken to act on functions that are
symmetric in the variables $z_i, z_j, z_k, \dots $ or $z_{j_1},z_{j_2}, \dots$
with
$j_{\kappa} \in J$, respectively.  When the restriction will be taken over all
the $N$
variables of the system, we shall simply use $\Res X$.
{}From the identities (3.3), it is checked that the
operators $D_i$ satisfy the commutation relations
\begin{equation}
[D_i,D_j]
	=\beta (D_j-D_i)K_{ij} .
\end{equation}
With $m$ some integer, it is then straightforward to verify that
\begin{multline}
\Res  ^{\{i,j\}} (D_i+m\beta )  (D_j+(m+1)\beta ) \\
 = \Res  ^{\{i,j\}}  (D_j+m\beta )  \bigl(D_i+(m+1)\beta \bigr).
\end{multline}
It follows that $\Res  ^J D_{J, k}= \Res  ^J  (D_{j_1}+k\beta ) \dots
\bigl(D_{j_{\ell}}+(k+\ell-1)\beta\bigr)$ is invariant under the permutations
of the
variables $z_{j_{\kappa}}$, $j_{\kappa} \in J$ and that this operator therefore
leaves
invariant the space of symmetric functions in these variables.  Recalling how
the $B_i^{+}$ are constructed in (3.9) in  terms of the operators $D_{1,J}$,
it is clear that $\phi_{\lambda}= ( B_{N-1}^{+} )^{\lambda_{N-1}} \dots
(B_{2}^{+})^{\lambda_2-\lambda_3} ( B_{1}^{+}) ^{\lambda_{1}-\lambda_2} \cdot
1$ is a symmetric function of the variables
$z_1, \dots ,z_N $.  That it is a homogeneous polynomial of degree $n$ is
readily
seen by observing that the operators $D_{J,k}$ have scaling dimension zero
and hence that $B_i^{+} \to \rho^i B_i^{+}$ when $z_i \to \rho z_i$,
$i=1, \dots ,N$.  The degree of $\phi_{\lambda}$ is thus
$\lambda_1 -\lambda_2 +2(\lambda_2-\lambda_3)+ \dots +(N-1)\lambda_{N-1}=n$.

\subsubsection{} It might be useful to give an example.  To this end, let us
take the
number of particles $N=3$ and consider a solution of degree $n=4$.  According
to  formulas (3.13)--(3.15), the Jack polynomial associated to the partition
$\lambda=(3,1)$ is given by
\begin{equation}
J_{(3,1)} \left( z_1,z_2,z_3;1/\beta \right)
	= c_{(3,1)}^{-1}  B_2^{+} (B_1^{+})^2 \cdot 1,
\end{equation}
with
\begin{equation}
\begin{split}
B_1^{+}
	&= z_1(D_1+\beta)+z_2(D_2+\beta)+z_3(D_3+\beta) ,\\
B_2^{+}
	&= z_1 z_2 (D_1+\beta)(D_2+2 \beta)+z_1 z_3 (D_1+\beta)(D_3+2 \beta)\\
	& \qquad\qquad\qquad\qquad + z_2 z_3 (D_2+\beta)(D_3+2 \beta),
\end{split}
\end{equation}
and
\begin{equation}
c_{(3,1)}
	=2\beta^2(\beta+1)^2.
\end{equation}

\subsubsection{}  We shall conclude this section by describing yet another way
to
obtain the spectrum and to characterize the wave functions of the CS model
with the help of operators closely related to the operators $D_i$ introduced
in (3.4).  Let
\begin{equation}
\widehat D_i
	= D_i + \beta (i-1) - \beta \sum_{j<i} (1-K_{ij}), \qquad i=1, \dots ,N.
\end{equation}
Remarkably, this defines a set of $N$ commuting operators
\begin{equation}
[\widehat D_i, \widehat D_j ]
	= 0.
\end{equation}
We also have
\begin{equation}
\begin{split}
[K_{ii+1},\widehat D_k]= & 0 \quad \text{ if } k \neq i, i+1\\
K_{ii+1} \widehat D_i- & \widehat D_{i+1} K_{ii+1}=\beta .
\end{split}
\end{equation}
In terms of these the Hamiltonian $H$ of (2.9) reads
\begin{equation}
H
	= \Res   \widehat H,
\end{equation}
with
\begin{equation}
\widehat H
	= \sum_{i=1}^N \Bigl\{ \widehat D_i^2 -(N-1) \beta \widehat D_i \Bigr\}+
\frac{1}{6} N(N-1)
(N-2) \beta^2.
\end{equation}
Clearly, $[\widehat H, \widehat D_i]=0 $ for all $i$.  We may thus obtain the
eigenfunctions of $\widehat H$ by diagonalizing simultaneously all the
$\widehat D_i$.
The eigenvalues and eigenfunctions of these operators are simply constructed
by observing that they are triangular on the set of monomials $\hat m_{\lambda}
= z_i^{\lambda_1} z_2^{\lambda_2} \dots z_N^{\lambda_N}$ associated to the
partitions $\lambda=(\lambda_1,  \dots ,\lambda_N)$ of $n$.  (Note that the
monomials $\hat m_{\lambda}$ are not symmetrized.)  The functions
$\chi_{\lambda}$ such that
\begin{equation}
\widehat D_i ~ \chi_{\lambda}
	= \delta_i^{\lambda} ~ \chi_{\lambda},
\end{equation}
have the form
\begin{equation}
\chi_{\lambda}
	=\hat m_{\lambda} + \sum_{\mu < \lambda} u_{\lambda \mu} \hat m_{\mu}, \qquad
|\lambda|=|\mu|=n,
\end{equation}
and from the action of $\widehat D_i $ on $\hat m_{\lambda}$, the eigenvalues
$\delta_i^{\lambda}$ are found to be
\begin{equation}
\delta_i^{\lambda}
	=\lambda_i + \beta(N-i).
\end{equation}
The eigenvalues of $\widehat H$ are then readily evaluated from (3.27) and one
gets
\begin{equation}
\widehat H  \chi_{\lambda}
	= \sum_{j=1}^N \bigl( \lambda_j^2+\beta(N+1-2j) \lambda_j
\bigr)  \chi_{\lambda}.
\end{equation}
Comparing with (2.22), we see that $H$ and $\widehat H$ have the same spectrum.

Using the relations (3.25), we show that
\begin{equation}
K_{ii+1}
	\widehat H= \widehat H  K_{ii+1}.
\end{equation}
Define now
\begin{equation}
\phi_{\lambda}
	=\sum_{\text{permutations}}\chi_{\lambda}.
\end{equation}
Since all permutations can be expressed as products of transpositions, owing to
(3.32), we see that
\begin{equation}
\widehat H  \phi_{\lambda}
	= \sum_{\text{permutations}} \widehat H  \chi_{\lambda}.
\end{equation}
We thus find that the symmetric function $\phi_{\lambda}$ which has triangular
expansion on the monomial basis $\bigl\{ m_{\lambda};|\lambda|=n \bigr\}$ is
actually an
eigenfunction of $\widehat H$ and also, of $H= \Res   \widehat H$ obviously,
with
eigenvalue equal to $\sum_{j=1}^N \bigl( \lambda_j^2+\beta(N+1-2j) \lambda_j
\bigr)$.  We thus conclude that it must be proportional to the Jack polynomial
\begin{equation}
J_{\lambda} \left( z_1, \dots  ,z_N;1/\beta \right) \sim
\sum_{\text{permutations}}
\chi_{\lambda} (z_1, \dots ,z_N;1/\beta).
\end{equation}
This provides an alternative to the method described in section~2 for obtaining
the spectrum and eigenfunctions of the CS Hamiltonian.

\section{Proofs}
\subsection{Outline}
After preliminary remarks, we shall proceed to give the proof of Theorem~3.1
which we stated in the last section.  To this end, we shall require additional
definitions.  Keeping with the notation of section~3, we introduce
\begin{equation}
N_{i,J}
	= \sum_{\begin{subarray}{c} J' \subset J\\ |J'| = i \end{subarray}} D_{0,J'},
\qquad i = 1, \dots , N,
\end{equation}
and the short hands
\begin{equation}
\begin{split}
&N_i
	\equiv N_{i, \{1, \dots ,N \}},\\
&\tilde N_{i,J}
	\equiv N_{i,\{ J_1,\dots , J_i \}}=D_{0,J},\\
&\tilde N_{i}
	\equiv N_{i,\{ 1,\dots , i \}}=D_{0,\{ 1,\dots,i \}}.
\end{split}
\end{equation}
We remind the reader that the notation $\Res^J X$ indicates that the operator
is taken to act on functions that are symmetric under the exchanges of the
variables $z_{j_{\kappa}}$, $j_{\kappa} \in J$.  At times, we shall need to
vary the set of variables entering in the Dunkl operators $\nabla_i$ from which
the operators $B_{i,J}^+$ and $N_{i,J}$ are built.  To specify this, if $S$ is
the set of integers labelling the variables, we shall introduce a superscript
$S$ writing for instance
\begin{equation}
\nabla_i^{(S)}
	= \frac{\partial}{\partial z_i} + \beta \sum_{\begin{subarray}{c} j \in S\\
j \ne i \end{subarray}} \frac{1}{z_i-z_j}   (1 - K_{ij}) ,\qquad i \in S.
\end{equation}
The same superscript will be added, when necessary, to the symbols of
the operators constructed from these $\nabla_i^{(S)}$.  When $S$ consists of
the first $M$ integers, $S = \{ 1,2, \dots , M \}$ we shall write
$N_{i,J}^{(S)} = N_{i,J}^{\{M\}}$.  Generally, when this superscript is
omitted, it is understood that the operators depend on the original $N$
variables $z_1, z_2, \dots ,z_N$.  The only exception will occur in
subsection~4.2 where we shall need to use $M \geq N$ variables and shall also
drop the superscript at some point.

By the reasoning of 3.3.1, it is clear that the operators $\Res^JB_{i,J}^+$
 and $\Res^J N_{i,J}$ are invariant under the permutations of the variables
$z_{j_{\kappa}}$, $j_{\kappa} \in J$.  We may also remark that
\begin{equation}
[ \Res^J N_{i, J}, \Res^J N_{k,J}] = 0.
\end{equation}
This is seen as follows.  The operators $\Res^J N_{i,J}^{(J)}$ are completely
symmetric under the permutations of the indices of the operators $D_i^{(J)}$.
They must therefore be combinations of the invariants $L_i^{(J)} = \sum_{k \in
J} (D_k^{(J)})^i$, $[L_i^{(J)},L_j^{(J)}] = 0$ and as a result, must also
commute among themselves,
$[\Res^J  N_{i,J}^{(J)},  \Res^J N_{k,J}^{(J)}] = 0$.  Since $[\Res^J  N_{i,J},
 \Res^J N_{k,J}]$ = $ \Res^J [ N_{i,J},
  N_{k,J}]$ and $[\Res^J  N_{i,J}^{(J)},  \Res^J N_{k,J}^{(J)}]$ =  $\Res^J
[N_{i,J}^{(J)},   N_{k,J}^{(J)}]$ have the same form in terms of Dunkl
operators  and since the commutation relations between these operators are not
affected by the number of variables, (4.4) must then be true.

Let
\begin{equation}
\varphi_{(\lambda_1,\dots , \lambda_i, 0, \dots )}
	\equiv (B_i^+)^{\lambda_i} \dots  (B_1^+)^{\lambda_1 - \lambda_2} \varphi_0,
\end{equation}
with $\varphi_0 = 1$.  We shall show that $J_{(\lambda_1,\dots ,
\lambda_{N-1})}(z_i;1/\beta) = c_\lambda^{-1} \varphi_{(\lambda_1,\dots ,
\lambda_{N-1})}$ by proving that the functions $\varphi_{(\lambda_1,\dots ,
\lambda_{N-1})}$ thus constructed are simultaneous eigenfunctions of the
Hamiltonian $H$ given in (2.8)--(2.10) and of the momentum operator $P=
\frac{2\pi}{L} \sum_j z_j \partial/\partial z_j$. This is in essence the
content of Theorem~3.1. There will remain to prove that the constant
$c_\lambda$ is indeed given by~(3.14) and  (3.15).  This will be done in
subsection~4.4. We have already observed in subsection~3.3.1 that
$\varphi_\lambda \equiv \varphi_{(\lambda_1,\dots , \lambda_{N-1})}$ is
$S_N$--invariant and that it is a homogeneous polynomial of degree $|\lambda| =
n$ in the variables $z_1, \dots , z_N \colon P\varphi_\lambda = \frac{2\pi n
}{L} \varphi_\lambda$.  Since $\varphi_\lambda$ is symmetric, this last
property can be expressed as follows
\begin{equation}
\left( \sum_{j=1}^N D_j \right) \varphi_{(\lambda_1,\dots , \lambda_{N-1})} =
\left( \sum_{j=1}^{N-1} \lambda_j \right) \varphi_{(\lambda_1,\dots ,
\lambda_{N-1})}.
\end{equation}
Theorem~3.1 will be proved by establishing the following proposition.

\begin{proposition}
\begin{multline}
[H, B_i^+] \varphi_{(\lambda_1,\dots , \lambda_i,0,\dots)}\\
= B_i^+ \left\{
2 \sum_{j=1}^{N-1} \lambda_j + i + \beta i(N-i) \right\}
\varphi_{(\lambda_1,\dots , \lambda_i,0,\dots)}, \\   \forall i \in \{ 1, \dots
, N-1 \}.
\end{multline}
\end{proposition}
If we assume that this proposition is true, it is then readily seen that the
successive applications of the operators $B_k^+$, $k = 1, \dots ,N-1$, on
$\varphi_0$ build the eigenfunctions $\varphi_{(\lambda_1, \lambda_2,\dots)}$
of $H$.  Let $\varepsilon_{(\lambda_1, \lambda_2,\dots)}$ be their eigenvalues.
 Since
\begin{equation}
B_i^+ \varphi_{(\lambda_1, \dots  , \lambda_i, 0, \dots)}
	= \varphi_{(\lambda_1+1, \dots  , \lambda_i+1, 0, \dots)} ,
\end{equation}
by definition, (4.7) is equivalent to
\begin{multline}
H  \varphi_{(\lambda_1+1, \dots  , \lambda_i+1, 0 \dots)}\\
	= \left\{ \varepsilon_{(\lambda_1, \dots  , \lambda_i, 0 ,\dots)} + 2
\sum_{j=1}^{N-1} \lambda_j + i + \beta i(N-i) \right\}
\varphi_{(\lambda_1+1, \dots  , \lambda_i+1, 0, \dots)}.
\end{multline}
Iterating with this formula, starting from the function $\varphi_0 = 1$ for
which $H\varphi_0 = 0$, one recursively finds that
\begin{equation}
H\varphi_{(\lambda_1, \dots , \lambda_{N-1})}
	= \left\{ \sum_{j=1}^{N-1} \Bigl( \lambda_j^2 + \beta(N+1 - 2j)\lambda_j
\Bigr) \right\} \varphi_{(\lambda_1, \dots  , \lambda_{N-1})} ,
\end{equation}
which is the desired result.

Instead of proving Proposition~4.1 directly, we observe that it follows from
the Propositions~4.2 and 4.3 given below.

\begin{proposition}
\begin{equation}
\tilde N_{i+1,J}~ \varphi_{(\lambda_1 , \dots , \lambda_i, 0, \dots)} = 0,
\qquad 1 \le i \le N-1.
\end{equation}
\end{proposition}

\begin{proposition}
\begin{multline}
[H,B_i^+] \varphi_{(\lambda_1 , \dots , \lambda_i, 0, \dots)}\\
	= \left\{ B_i^+ \Bigl[ 2 \sum_{j=1}^N D_j  + i + \beta i(N-i) \Bigr] +
\sum_J G_J \tilde N_{i+1,J} \right\} \varphi_{(\lambda_1 , \dots , \lambda_i,
0, \dots)}, \\
 1 \le i \le N-1,
\end{multline}
where $G_J$ are certain unspecified expressions.
\end{proposition}

It is indeed  seen with the help of (4.6), that Proposition~4.3 is equivalent
to
Proposition~4.1 if Proposition~4.2 is true.  The next two subsections are
devoted
to the proofs of these Propositions~(4.2  and 4.3) which imply Proposition~4.1.
 We proceed by induction mostly; the symmetry properties of our constructs are
extensively used and formula (3.19) is for instance called upon repeatedly.
The constants $c_\lambda$ which relate the functions $\varphi_\lambda$ to the
monic Jack polynomials $J_\lambda$ are determined in the last subsection to
complete the proof of Theorem~3.1.

\subsection{Proof of Proposition~4.2}
Owing to the fact that $\varphi_{(\lambda_1 , \dots , \lambda_i, 0, \dots)}$ is
a symmetric function of the variables $z_1, \dots ,z_N$, it suffices to show
that
\begin{equation}
\tilde N_{ i+1} \varphi_{(\lambda_1 , \dots , \lambda_i, 0, \dots)} = 0,
\end{equation}
in order to prove (4.11).  Equation~(4.13) in turn, is seen to follow from the
relation
\begin{equation}
\Res[\tilde N_{i+1}, B_k^+] \sim \Res \tilde N_{k+1}, \qquad \forall k < i+1,
\end{equation}
{\sl where by $\sim$ we mean that the term on the right-hand side can be
multiplied
on the left by some non-singular operator}.  In fact, using (4.14) iteratively,
one
finds that $\tilde N_{i+1} \varphi_{(\lambda_1 , \dots , \lambda_i, 0, \dots)}
\sim \tilde N_{j+1} \varphi_0$ with $j$ the smallest integer such that
$\lambda_j - \lambda_{j+1} \ne 0$.  Since $\tilde N_{ j+1}\varphi_0 = 0$,
$\forall j$, (4.13) is thus implied by (4.14).

With the help of (3.19), it is easy to show that
\begin{equation}
\Res^{\{i+1\}} \tilde N_{i+1}
	= \Res^{\{i+1\}} \left( D_1 K_{1i+1} + \beta \sum_{j=2}^i K_{ji+1} +
\beta \right)\tilde N_{i},
\end{equation}
and hence that
\begin{equation}
\Res^{\{i+1\}} \tilde N_{i+1} \sim \Res^{\{i+1\}} \tilde N_{i}.
\end{equation}
Equation (4.14) will therefore be true if
\begin{equation}
\Res[\tilde N_{k+1}, B_k^+] \sim \Res \tilde N_{k+1},
\end{equation}
holds.  We shall actually prove a stronger result, namely:

\begin{proposition}
\begin{equation}
\Res^{\{N\}} [\tilde N_{k+1}^{\{M\}}, B_k^{+\{M\}} ] \sim \Res^{\{N\}}
\tilde N_{k+1}^{\{M\}}, \qquad \forall M \ge N.
\end{equation}
\end{proposition}

As already mentioned, the superscript $\{ M \}$ indicates that we are using
Dunkl operators that not only depend on the variables $z_1, \dots , z_N$ but
also on the variables $z_{N+1}, \dots , z_M$.  (It should be noted that $B_k^{+
\{M\}}$ is symmetric under $S_N$ but not under $S_M$.)  Clearly, proving
Proposition~4.4 is tantamount to proving Proposition~4.2.

For the remainder of subsection~4.2, the Dunkl operators entering our various
expressions will always be taken to depend on the variables $z_1, \dots ,z_M$.
With this understood, we shall omit the superscript $\{M \}$ in the following.

It is possible to isolate the parts of $B_k^+$ involving the indices $\{ 1,
\dots , k+1\}$.  From the definition of $B_k^+$, one has the result:

\begin{lemma}
\begin{equation}
B_k^+
	= \sum_{\ell = 0}^k \sum_{\begin{subarray}{l} J\subset \{k+2, \dots , N\}
\\ ~~~~|J| = \ell \end{subarray}} z_J B^+_{k-\ell,\{1, \dots , k+1 \}}
D_{k-\ell+1,J},
\end{equation}
with  $B_{0,J}^+ \equiv 1$ and when $|J|=0$, $z_J\equiv 1$ and
$D_{k+1,J}\equiv 1$.
\end{lemma}

We shall also make use of the following formulas.

\begin{lemma}
\begin{subequations}
\begin{align}
\mathrm{(i)} \ &[D_i, z_J]
	= -\beta z_i \sum_{j\in J} z_{J\backslash \{j\}} K_{ij}, \qquad i \notin J ,\\
\mathrm{(ii)} \ &[D_i, z_J]
	= z_J \biggl( 1 + \beta \sum_{\begin{subarray}{l} ~~~~j \notin J \\ j \in
\{1,\dots , M\} \end{subarray}} K_{ij} \biggr), \qquad i \in J, \\
\mathrm{(iii)} \ & \Res^{\{k+1\}} [\tilde N_{k+1}, z_\ell]\\
&= \Res^{\{k+1\}} (-\beta)(z_1K_{1 2}K_{23} \dots K_{k k+1}K_{k+1 \ell} +
 \dots + z_{k+1} K_{k+1 \ell})\tilde N_{k} \notag\\
&\qquad\qquad\qquad\qquad \qquad \sim \Res^{\{k+1\}} \tilde N_{k}, \qquad \ell
\notin \{ 1, \dots , k+1 \}. \notag
\end{align}
\end{subequations}
\end{lemma}
Equations (4.20a) and (4.20b) are immediately obtained from the definition of
the operators $D_i$.  Result (i) expresses the fact that the commutator of
$D_i$ with products of the form $z_{j_1}z_{j_2} \dots  z_{j_\ell}$ that exclude
$z_i$, yields expressions which have as factors on the left, similar products
with
 $z_i$ replacing one of the initial variables.  Formula (ii) shows that the
commutator of $D_i$ with products of the form $z_{j_1} z_{j_2} \dots
z_{j_\ell}$
involving $z_i$, gives expressions having as factor on the left the same
products
$z_{j_1} z_{j_2} \dots  z_{j_\ell}$ of the variables.

Property (iii) is easily derived by induction.  It shows that on symmetric
functions
of $z_1, \dots , z_{k+1}$, commuting $\tilde N_{k+1}$ with a variable $z_\ell$
not
belonging to this set, has the effect of giving an operator having $\tilde
N_{k}$
as factor on the right.

{\sl In the following, we shall need to identify in various expressions, the
terms that do not have $z_1$ appearing as an explicit factor on the left.  If
$X$ represents one such quantity of interest, the terms in question will be
denoted by $X\bigm|_{z_1 \sim 0}$}.  A result that will soon be useful is:

\begin{lemma}
\begin{multline}
\Bigl( [\tilde N_{ \ell+1},z_2 \dots z_n]
	+ z_2 \dots z_n [D_1, D_{1,\{2,\dots ,\ell+1\}}] \Bigr)\Bigm|_{z_1\sim 0} \sim
D_1, \\  \forall M \ge n \ge \ell+1.
\end{multline}
\end{lemma}

The proof is done by induction.  That (4.21) is true when $\ell = 1$ is easily
seen from
\begin{multline}
\Bigl( [\tilde N_{2}, z_2 \dots z_n] + z_2 \dots  z_n[D_1,D_2] \Bigr)
\Bigm|_{z_1 \sim 0}\\
	= z_2 \dots  z_n \left( 1 + \beta K_{12} + \beta \sum_{i = 3}^M K_{2i} \right)
D_1 \sim D_1.
\end{multline}
Now in general one has,
\begin{equation}
\begin{split}
&\Bigl( [\tilde N_{ \ell + 1}, z_2 \dots z_n] + z_2 \dots z_n[D_1,
D_{1,\{2,\dots , \ell + 1 \}}] \Bigr) \Bigm|_{z_1 \sim 0} \\
&~~~~~~= \biggl( [\tilde N_{\ell}, z_2 \dots z_n] (D_{\ell + 1} + \ell\beta) +
\tilde N_{ \ell}
 [D_{\ell + 1}, z_2  \dots  z_n]\\
&\qquad\quad\quad + z_2 \dots  z_n[D_1, D_{1,\{2,\dots , \ell \}}] (D_{\ell
+1}+ \ell \beta) \\
&\qquad\qquad\qquad\qquad\quad + z_2 \dots z_n D_{1,\{2, \dots ,\ell\}} [D_1,
D_{\ell+1}] \biggr) \biggm|_{z_1 \sim 0}.
 \end{split}
\end{equation}
It is easy to check that
\begin{equation}
\begin{split}
&\tilde N_{ \ell}[D_{\ell+1}, z_2 \dots  z_n]
	= z_2 \dots  z_n \tilde N_{ \ell} \biggl( 1 + \beta K_{1\ell + 1} +
\beta \sum_{i=n+1}^M K_{i \ell+1}\biggr)\\
&\quad\quad + [\tilde N_{\ell}, z_2 \dots  z_n]
\biggl( 1 + \beta K_{1 \ell + 1} + \beta \sum_{i=n+1}^M K_{i \ell+1}\biggr)\\
&=z_2 \dots  z_n \tilde N_{\ell}(\beta K_{1\ell+1}) + [\tilde N_{\ell}, z_2
\dots  z_n ]
\biggl( 1 + \beta K_{1\ell + 1} + \beta \sum_{i=n+1}^M K_{i \ell+1}\biggr)\\
&\quad\quad + z_2 \dots  z_n \bigl( D_{1,\{2,\dots , \ell \}} D_1 +
[D_1, D_{1,\{2,\dots , \ell\}}] \bigr)
\biggl( 1 + \beta  \sum_{i=n+1}^M K_{i \ell+1}\biggr)
\end{split}
\end{equation}
and that
\begin{multline}
D_{1,\{2, \dots , \ell \}} [D_1,D_{\ell+1}]
 	=D_{1,\{2, \dots , \ell \}} (\beta K_{1 \ell+1} D_1)\\
-\tilde N_{\ell} (\beta K_{1\ell+1}) + [D_1, D_{1,\{2, \dots , \ell \}}] (\beta
K_{1 \ell+1}).
\end{multline}
Substituting (4.24) and (4.25) in (4.23) and dropping terms which have $D_1$
already on the right, we get
\begin{equation}
\begin{split}
&\Bigl( [\tilde N_{ \ell + 1}, z_2 \dots z_n] + z_2 \dots  z_n[D_1,
D_{1,\{2,\dots , \ell + 1 \}}] \Bigr) \Bigm|_{z_1 \sim 0} \\
&\qquad\qquad \longrightarrow \biggl( \bigl\{ [\tilde N_{\ell}, z_2 \dots z_n]
+ z_2  \dots  z_n [D_1, D_{1,\{2,\dots , \ell \}}] \bigr\} \\
&\qquad\qquad\qquad\qquad \times \Bigl( D_{\ell+1} + \ell \beta + 1 +
\beta K_{1 \ell+1} + \beta \sum^M_{i=n+1} K_{i \ell+1} \Bigr) \biggr)
\biggm|_{z_1\sim 0}.
\end{split}
\end{equation}

By hypothesis, the term in curly brackets has $D_1$ occurring on the right and
since
\begin{equation}
D_1(D_{\ell + 1} + \beta K_{1 \ell + 1}) = (D_{\ell + 1} + \beta K_{1 \ell +
1})D_1,
\end{equation}
Lemma~4.7 is thus shown to hold.

As a step toward establishing Proposition ~4.4, we shall  prove the following
result.

\begin{proposition}
\begin{equation}
\Res^{\{k+1\}} [\tilde N_{k+1}, B^+_{k,\{1,\dots , k+1\}}] \sim \Res^{\{k+1\}}
\tilde N_{k+1}.
\end{equation}
\end{proposition}
\noindent Note that (4.28) is a special case of (4.18) with $N=k+1$.

Let,
\begin{equation}
z_1 \dots \hat z_i \dots z_{k+1}
	= \prod^{k+1}_{\begin{subarray}{l} \ell = 1\\ \ell \ne i \end{subarray}}
z_\ell.
\end{equation}
{}From the identities (i) and (ii) of Lemma~4.6, it is immediate to see that
\begin{equation}
[\tilde N_{k+1}, z_1 \dots  \hat z_i \dots  z_{k+1}]
	= \sum_{j=1}^{k+1} z_1 \dots  \hat z_j  \dots  z_{k+1} \mathcal E_j,
\end{equation}
with $\mathcal E_j$ quantities involving the operators $D_\ell$ and $K_{\ell
m}$
, $1 \le \ell, m \le k+1$.  It thus follows, given the definition of
 $B_{k,\{1,\dots , k+1\}}^+$, that all the terms in the expression of the
commutator
 $[\tilde N_{k+1}, B_{k,\{1,\dots , k+1\}}^+]$ will have on the left a factor
consisting in the product of $k$ distinct variables taken among the set
$\{z_1, z_2, \dots , z_{k+1}\}$.  Since $\tilde N_{k+1}$ and
$B_{k,\{1,\dots , k+1\}}^+$ are invariant under permutations of these
variables,
in order to prove Proposition~4.8, it will suffice to show that the term in the
expression of this commutator which is multiplied on the left by
$z_2 \dots  z_{k+1}$ has on the right, the operator $\Res^{\{k+1\}}
\tilde N_{k+1}$.  In other words, proving that
\begin{equation}
\Res^{\{k+1\}} [\tilde N_{k+1}, B_{k,\{1,\dots , k+1\}}^+] \bigm|_{z_1 \sim 0}
\sim
\Res^{\{k+1\}} \tilde N_{k+1},
\end{equation}
will establish Proposition~4.8.  We proceed by induction.  When $k=1$ we have,
\begin{equation}
\Res^{\{2\}} [\tilde N_{2}, B_{1,\{1,2\}}^+ ] \bigm|_{z_1 \sim 0}
	= \Res^{\{2\}} z_2 \left( 1 + \beta \sum_{j=3}^M K_{2 j} \right) \tilde N_{2}
\sim \Res^{\{2\}} \tilde N_{2}.
\end{equation}
We now suppose that (4.31) is true for all $k$ smaller than $\ell$ hoping  that
it is satisfied for $k = \ell$ as a consequence.  Let us cast $B_{\ell, \{ 1,
\dots , \ell+1\}}^+$ in a form where the part involving $z_1$ explicitly is
isolated:
\begin{equation}
B_{\ell, \{ 1, \dots , \ell+1\}}^+
	= z_1 B_{\ell-1, \{ 2, \dots , \ell+1\}}^+ (D_1+\ell \beta)
	+ z_2  \dots  z_{\ell+1} D_{1,\{2, \dots , \ell +1 \}}.
\end{equation}
With the help of (4.33) and since $D_1 \cdot D_{1,\{2, \dots , \ell +1\}} =
\tilde N_{ \ell +1}$,  we can write:
\begin{equation}
\begin{split}
&\Res^{\{\ell+1\}} [\tilde N_{ \ell+1}, B_{\ell,\{1, \dots , \ell+1\}}^+ ]
\bigm|_{z_1 \sim 0}\\
&~~~~~= \Res^{\{\ell + 1\}} \biggl( D_1 [D_{1,\{2, \dots , \ell+1\}}, z_1]
B_{\ell-1, \{ 2, \dots , \ell+1\}}^+ (D_1 + \ell \beta)\\
&\qquad + \Bigl\{ [\tilde N_{ \ell+1}, z_2 \dots  z_{\ell +1}]
 + z_2 \dots z_{\ell + 1} [D_1,D_{1,\{2, \dots , \ell+1\}}] \Bigr\}
 D_{1,\{2, \dots , \ell + 1\}}\biggr) \biggm|_{z_1\sim 0}.\\
& ~~~~~~~~~~
\end{split}
\end{equation}

{}From Lemma~4.7, we see that the term in curly brackets has
the operator $D_1$ as last factor on the right.
It thus remains to show that
\begin{equation}
\Res^{\{\ell +1\}} \biggl( D_1[D_{1,\{2, \dots ,\ell+1\}}, z_1]
B_{\ell - 1, \{2, \dots , \ell + 1\}}^+ (D_1 + \ell \beta) \biggr) \biggm|_{z_1
\sim 0}
 \sim \Res^{\{\ell +1\}} \tilde N_{\ell+1},
\end{equation}
to complete the proof of Proposition~4.8.  A result analogous to (4.20c) which
is again straightforwardly proved by induction is now helpful:
\begin{equation}
\begin{split}
&\Res^{\{2, \dots , \ell+1\}} [D_{1,\{2,\dots , \ell \}}, z_1]\\
&\quad \quad = \Res^{\{2, \dots , \ell+1\}}(-\beta)(z_2 K_{2 3} \dots  K_{\ell
\ell + 1}
K_{1 \ell + 1} + \dots + z_{\ell + 1} K_{1 \ell+1})
 D_{1,\{2,\dots , \ell \}}.\\
&
\end{split}
\end{equation}
Since $B_{\ell - 1, \{2, \dots , \ell + 1\}}^+ (D_1 + \ell \beta)$ is invariant
under the permutations of the indices $\{ 2, \dots ,\ell + 1 \}$, we find with
the help of (4.36) that
\begin{multline}
\Res^{\{ \ell + 1\}} \biggl( D_1 [D_{1,\{2,\dots , \ell + 1\}}, z_1]
B_{\ell - 1, \{2, \dots , \ell + 1\}}^+ (D_1 + \ell \beta) \biggr) \biggm|_{z_1
\sim 0}\\
\sim \Res^{\{ \ell + 1\}} \biggl( \tilde N_{\ell,\{2,\dots , \ell + 1\}}
B_{\ell - 1, \{2, \dots , \ell + 1\}}^+ (D_1 + \ell \beta) \biggr) \biggm|_{z_1
\sim 0},
\end{multline}
where in obtaining (4.37), we have used the identity
\begin{equation}
\begin{split}
&D_1 \Res^{\{ 2, \dots ,\ell + 1\}}
(z_2K_{2 3} \dots  K_{\ell \ell+1} K_{1 \ell+1} + \dots + z_{\ell+1} K_{1 \ell
+ 1}) D_{1,\{2,\dots,\ell\}} \bigm|_{z_1 \sim 0}\\
 &= \Res^{\{ 2, \dots ,\ell + 1\}}
(z_2K_{2 3} \dots  K_{\ell \ell+1} K_{1 \ell+1} + \dots + z_{\ell+1} K_{1 \ell
+ 1 }) D_{\ell +1} D_{1,\{2,\dots,\ell\}} \bigm|_{z_1 \sim 0},
\end{split}
\end{equation}
which follows from the fact that $\Res^{\{\ell + 1\}} \Res^{\{2,\dots ,
\ell+1\}}
= \Res^{\{\ell+1\}}$ and that every term with which $D_1$ is commuted contains
the operator $K_{1 \ell+1}$.

{}From the induction hypothesis we have
\begin{equation}
\Res^{\{2,\dots , \ell + 1\}}
[\tilde N_{\ell, \{2 , \dots , \ell +1\}}, B_{\ell - 1, \{ 2, \dots , \ell
+1\}}^+ ] \bigm|_{z_1 \sim 0}
\sim \Res^{\{2,\dots , \ell + 1\}}\tilde N_{\ell, \{2 , \dots , \ell +1\}}.
\end{equation}
Since $(D_1 + \ell \beta)$ is invariant under the permutations of the indices
$\{2, \dots , \ell +1 \}$ we finally obtain
\begin{multline}
\Res^{\{2,\dots , \ell + 1\}}
\biggl( \tilde N_{\ell, \{2 , \dots , \ell +1\}} B_{\ell - 1, \{ 2, \dots ,
\ell +1\}}^+
(D_1 + \ell \beta) \biggr) \biggm|_{z_1 \sim 0}\\
\sim \Res^{\{2,\dots , \ell + 1\}}\tilde N_{\ell, \{2 , \dots , \ell +1\}}
(D_1 + \ell\beta) \sim \Res^{\{2,\dots , \ell + 1\}}\tilde N_{  \ell +1} ,
\end{multline}
which through (4.37) proves (4.35) and, as a consequence, Proposition~4.8.

A few more results will be required to prove Proposition~4.4.  They can be
stated as follows.

\begin{lemma}
For sets $J = \{ j_1, \dots, j_\ell \}$ of cardinality $\ell$ such that $J
\cap$~$ \{ 1, \dots , k+1 \} = \varnothing$, we have
\begin{subequations}
\begin{align}
\mathrm{(i)} \ &\Res^{\{k+1\}} [ \tilde N_{k+1}, z_J] \sim  \Res^{\{k+1\}}
\tilde N_{k+1-\ell} ,\\
\mathrm{(ii)} \  &\Res^{\{N\}} \bigl( \Res^{\{k+1\}} \tilde N_{k+1-\ell}
D_{k+1-\ell,J} \bigr) \sim \Res^{\{N\}} \tilde N_{k+1}, \\
\mathrm{(iii)} \ & \Res^{\{N\}} [\tilde N_{k+1}, D_{k+1-\ell,J}] \sim
\Res^{\{N\}}  \tilde N_{k+1}.
\end{align}
\end{subequations}
\end{lemma}

Property (i) is readily obtained by applying (4.20c) $\ell$ times.  Property
(ii) is
also easily derived.  Since $D_{k+1-\ell,J}$ is invariant under the
permutations
 of $\{1,\dots , k+1\}$, the restriction $\Res^{\{k+1\}}$ is redundant and can
be
dropped.  By definition, $\tilde N_{ k+1-\ell} D_{k+1-\ell,J} =
\tilde N_{k+1, \{1,\dots , k+1-\ell \} \cup J}$.  When restricting this
operator to
symmetric functions of $z_1, \dots , z_N$, we can relabel appropriately the
variables to find $\Res^{\{N\}} \tilde N_{k+1}$.  The last property is obtained
as
follows.  Using the identity
\begin{equation}
\Res^{\{1,\dots ,k+1,i\}} [ \tilde N_{k+1}, D_i] \sim \Res^{\{1, \dots
,k+1,i\}}
\tilde N_{k+1}, \qquad i \notin \{ 1, \dots , k+1 \},
\end{equation}
which is easy to prove, one sees that
\begin{equation}
\begin{split}
&\Res^{\{N\}}  \tilde N_{k+1} D_{k+1-\ell, J} \\
&\quad ~~~ = \Res^{\{N\}} \Res^{\{1,\dots , k+1, J_1\}}  \biggl( \tilde N_{k+1}
\bigl( D_{J_1} +
(k+1-\ell)\beta \bigr) \biggr) D_{k+2-\ell, J\backslash \{J_1\}}\\
&\qquad\qquad \qquad \sim \Res^{\{N\}} \tilde N_{k+1} D_{k+2-\ell, J\backslash
\{J_1\}}.
\end{split}
\end{equation}
Iterating (4.43), one thus arrives at (4.41c).

We are now ready to give the proof of Proposition~4.4.  We proceed again by
induction.  The proposition is easily seen to hold when $k=1$.  Indeed,
\begin{multline}
\Res^{\{N\}} [\tilde N_{2},B_1^+]
	= \Res^{\{N\}} [\tilde N_{2}, B_{1,\{1,2\}}^+] \\
	+ \sum_{i=3}^N \Res^{\{N\}} \bigl[ \tilde N_{2}, z_i(D_i+\beta)
\bigr].~~~~~~~~~~~~~~
\end{multline}

{}From Proposition~4.8 we know that $\Res^{\{N\}} [\tilde N_{2},
B_{1,\{1,2\}}^+]
\sim \Res^{\{N\}} \tilde N_{2}$.  Using (4.20c), we find
\begin{equation}
 \Res^{\{N\}} [\tilde N_{2}, z_i] (D_i + \beta) \sim \Res^{\{N\}} D_1 (D_i +
\beta)
\sim \Res^{\{N\}} \tilde N_{2},
\end{equation}
and, with the help of (4.41c), we observe that
\begin{equation}
\Res^{\{N\}} [\tilde N_{2}, D_i] \sim \Res^{\{N\}} \tilde N_{2}.
\end{equation}
All terms are therefore seen to have the required factor $\tilde N_{2}$ on the
right.

Now, upon supposing that $\Res^{\{N\}} [\tilde N_{k+1}, B_k^+] \sim
\Res^{\{N\}}
\tilde N_{k+1}$ is true for all $k < m$, we wish to prove that this relation
holds
also for $k = m$.  In view of  formula (4.19) for $B_m^+$, it is clear that
Proposition~4.4 would be established if  one could show for all $\ell \in \{0,
\dots , m\}$ that
\begin{equation}
\Res^{\{N\}} [\tilde N_{m+1}, z_J B_{m-\ell, \{1,\dots ,m+1\}}^+ D_{m-\ell+1,J}
]
\sim \Res^{\{N\}} \tilde N_{m+1},
\end{equation}
where $J$ are subsets of $\{m+2, \dots , N\}$ with cardinality $\ell$.  In
fact, property (4.41c) shows that it is sufficient to prove that
\begin{equation}
\Res^{\{N\}} [\tilde N_{m+1}, z_J B_{m-\ell, \{1,\dots ,m+1\}}^+]
D_{m-\ell+1,J}
\sim \Res^{\{N\}} \tilde N_{m+1},
\end{equation}
in order to establish (4.47).

The case $\ell = 0$ is the content of Proposition~4.8 and has thus already been
proved.  Remarkably, the cases of lower degree follow from the induction
hypothesis which allows one to assume that
\begin{equation}
\Res^{\{m+1\}} [\tilde N_{m-\ell+1},  B_{m-\ell, \{1,\dots ,m+1\}}^+]
 \sim \Res^{\{m+1\}} \tilde N_{m-\ell+1},
\end{equation}
for $\ell = 1,2, \dots , m$.

Owing to the invariance of $B_{m-\ell, \{1,\dots ,m+1\}}^+ D_{m-\ell+1,J}$ and
of  $D_{m-\ell+1,J}$ under the permutations of the indices $\{ 1, \dots , m+1
\}$, we may write
\begin{equation}
\begin{split}
&\Res^{\{N\}} [ \tilde N_{m+1}, z_J  B_{m-\ell, \{1,\dots
,m+1\}}^+]D_{m-\ell+1,J}\\
& = \Res^{\{N\}} \biggl( \Res^{\{m+1\}} [ \tilde N_{m+1}, z_J] \biggr)
B_{m-\ell, \{1,\dots ,m+1\}}^+ D_{m-\ell+1,J}\\
&\qquad\qquad + \Res^{\{N\}} z_J  \Bigl( \Res^{\{m+1\}} [ \tilde N_{m+1},
B_{m-\ell, \{1,\dots ,m+1\}}^+ ] \Bigr)D_{m-\ell+1,J}.
\end{split}
\end{equation}
We then use (4.41a), (4.49) and (4.41b) to conclude that (4.48) is true and
thus to
 finally complete the proof of Propostion~4.4.

The reason why we needed to prove the stronger result (4.18) instead of the
weaker one (4.17) should now have become clear.  Indeed, the relations (4.49)
with precisely the restriction $\Res^{\{m+1\}}$, are required for the induction
 proof of Proposition~4.4 to work.  These relations obviously follow from
 assuming that $\Res^{\{N\}} [\tilde N_{k+1}^{\{M\}}, B_{k,\{1,\dots,N\}}^{+
\{M\}}] \sim
\Res^{\{N\}}
 \tilde N_{k+1}^{\{M\}}$, $\forall k < m$ and $\forall M \ge N$.  One simply
takes
$N = m+1$ and $k = m-\ell$, $\ell = 1,2, \dots ,m$.  They would not have been
legitimate assumptions however, had the induction been performed on the
relation
$\Res^{\{N\}} [\tilde N_{k+1}^{\{N\}}, B_{k,\{1,\dots,N\}}^{+\{N\}}] \sim
\Res^{\{N\}}
\tilde N^{\{N\}}_{k+1}$, since in this case, the restriction is tied to the
number of
 variables $N$.

\subsection{Proof of Proposition~4.3}

The proof of Proposition 4.3 will be done in two steps.  We shall first show
that (4.11) holds when $i=N-1$ and the number of variables $ z_1,z_2,\dots$ is
equal to $N$.  We shall then establish Proposition 4.3 in full generality by
demonstrating that it is also true when the number of variables is taken to be
arbitrarily larger than $N-1$.  We shall need the following results to proceed.

\begin{lemma}
Let $z_1 \dots \hat z_i \dots z_{N}
	= \prod^{N}_{\begin{subarray}{l} \ell = 1\\ \ell \ne i \end{subarray}}
z_\ell$.  The following relations are satisfied.
\begin{subequations}
\begin{align}
\mathrm{(i)} \ &[D_j^2,z_1 \dots \hat z_i \dots z_{N} ] = \notag\\
&\qquad z_1 \dots \hat z_i \dots z_{N} \Bigl\{ (1+\beta K_{ij})^2 +(1+\beta
K_{ij})D_j+D_j(1+\beta K_{ij}) \Bigr\} ~~~~i \neq j \\
\mathrm{(ii)} \ &[D_i^2,z_1 \dots \hat z_i \dots z_{N} ] \bigm|_{z_1 \sim
0}=\notag\\
&\qquad  -\beta z_2 \dots z_N ( D_i K_{1i}+K_{1i} D_i + \beta K_{1i}+1)
\end{align}
\end{subequations}
\end{lemma}
Both identities are obtained straightforwardly from (4.20a) and (4.20b).
Let us introduce the notation
\begin{equation}
H_M = \sum_{i=1}^M D_i^2.
\end{equation}
The Hamiltonian in $N$ dimension $H^{\{N\}}$ (or $H$) is $H^{\{N\}} \equiv
\Res H_N^{\{N\}}$.  The following result is an immediate consequence of Lemma
4.10.
\begin{corollary}  For all $i$ strictly smaller than $N$, the operators  $
[~H_{N-1}^{\{N\}}$ , $z_1 \dots \hat z_i \dots z_{N}] \bigm|_{z_1 \sim 0}$
and $z_N [H_{N-1}^{\{N-1\}},z_1 \dots \hat z_i \dots z_{N-1}] \bigm|_{z_1 \sim
0}$
have the same symbol structure, that is, they have the same form in terms of
the coordinates and the operators $D_i$.
\end{corollary}
\begin{lemma}
\begin{subequations}
\begin{align}
\mathrm{(i)} \ & \Res^{\{N\}} [\sum_{i=1}^{N-1} D_i^{\ell},D_N]=\Res^{\{N\}}
\beta \Bigl( (N-1)D_N^{\ell}-\sum_{i=1}^{N-1} D_i^{\ell} \Bigr),\\
\mathrm{(ii)} \ & \Res^{\{N\}} [ D_N^{\ell},D_{1,\{1,\dots,N\}}]=  \\
&\qquad \Res^{\{N\}} \beta D_{1,\{1,\dots,N-1\}} \Bigl(
\sum_{i=1}^{N-1}D_i^{\ell}-(N-1)D_N^{\ell} \Bigr) ~~~\forall N \geq 2.  \notag
\end{align}
\end{subequations}
\end{lemma}
Formula (4.53a) is readily obtained from observing that
\begin{equation}
[D_i^{\ell},D_j]=\beta (D_j^{\ell}-D_i^{\ell})K_{ij}.
\end{equation}
The proof of (4.53b) is more involved and proceeds by induction.  With the help
of (4.54),(4.53b) is easily seen to hold in the first non-trivial case:
\begin{equation}
\begin{split}
\Res^{2} [ D_2^{\ell},(D_1+\beta)(D_2+2\beta)]&=\Res^{\{2\}}\beta
(D_1^{\ell}-D_2^{\ell})(D_1+2\beta)\\
& =\Res^{\{2\}}\beta (D_1+\beta)(D_1^{\ell}-D_2^{\ell}).
\end{split}
\end{equation}
We show that (4.53b) follows from assuming that
\begin{equation}
\begin{split}
& \Res^{\{1,\dots,N-2,N\}}[D_N^{\ell},D_{1,\{1,\dots,N-2,N\}}]=\\
& \qquad \qquad \Res^{\{1,\dots,N-2,N\}} \beta D_{1,\{1,\dots,N-2\}}\Bigl(
\sum_{i=1}^{N-2} D_i^{\ell} - (N-2)D_N^{\ell} \Bigr)
\end{split}
\end{equation}
is valid $\forall~N>3$.

Simple manipulations give
\begin{equation}
\begin{split}
& \Res^{\{N\}}[D_N^{\ell},D_{1,\{1,\dots,N\}}]=\\
& \qquad \qquad
\Res^{\{N\}}[D_N^{\ell},D_{1,\{1,\dots,N-2,N\}}(D_{N-1}+N\beta)]\\
& \qquad \qquad \qquad \qquad- \beta
\Res^{\{N\}}D_{1,\{1,\dots,N-2\}}(D_{N-1}-D_N)(D_N^{\ell}-D_{N-1}^{\ell}).
\end{split}
\end{equation}
One then easily finds that
\begin{equation}
\begin{split}
& \Res^{\{N\}}[D_N^{\ell},D_{1,\{1,\dots,N-2,N\}}(D_{N-1}+N \beta)]\\
& \qquad \qquad =\Res^{\{N\}}\Bigl( \Res^{\{1,\dots,N-2,N\}}
[D_N^{\ell},D_{1,\{1,\dots,N-2,N\}}]\Bigr)(D_{N-1}+N \beta)\\
& \qquad \qquad \qquad \qquad \qquad \qquad
+\Res^{\{N\}}D_{1,\{1,\dots,N-2,N\}} [D_N^{\ell},D_{N-1}]\\
& \qquad \qquad = \beta \Res^{\{N\}}D_{1,\{1,\dots,N-2\}} \biggl\{ \Bigl(
\sum_{i=1}^{N-2} D_i^{\ell} - (N-2)D_N^{\ell}\Bigr) (D_{N-1}+N \beta)\\
& \qquad \qquad \qquad \qquad \qquad \qquad+ \bigl( D_N + (N-1)\beta \bigr)
(D_{N-1}^{\ell}-D_N^{\ell}) \biggr\}.
\end{split}
\end{equation}
The induction hypothesis (4.56) has been used in (4.58) to obtain the last
equality.  From (4.53a) we also get
\begin{equation}
\begin{split}
& \Res \Bigl( \sum_{i=1}^{N-2} D_i^{\ell} - (N-2)D_N^{\ell}\Bigr) (D_{N-1}+N
\beta)=\\
& \qquad \qquad \qquad \qquad \qquad \Res \bigl( D_{N-1}+(N-1)\beta \bigr)
\Bigl( \sum_{i=1}^{N-2} D_i^{\ell}
- (N-2)D_N^{\ell}\Bigr).
\end{split}
\end{equation}
Inserting this last relation in (4.58) and reorganizing the terms yields
\begin{equation}
\begin{split}
& \Res^{\{N\}}[D_N^{\ell},D_{1,\{1,\dots,N-2,N\}}(D_{N-1}+N \beta)]=\\
& \qquad \qquad ~~~~ ~~~~\beta \Res^{\{N\}} D_{1,\{1,\dots,N-1\}}\Bigl(
\sum_{i=1}^{N-1} D_i^{\ell} - (N-1)D_N^{\ell}\Bigr)\\
& \qquad \qquad ~~~~-\beta \Res^{\{N\}} D_{1,\{1,\dots,N-2\}} \bigl(
D_{N-1}+(N-1) \beta \bigr)(D_{N-1}^{\ell}-D_N^{\ell}) \\
& \qquad \qquad ~~~~+\beta  \Res^{\{N\}} D_{1,\{1,\dots,N-2\}} \bigl(
D_{N}+(N-1) \beta \bigr)(D_{N-1}^{\ell}-D_N^{\ell}) .
\end{split}
\end{equation}
Using this result in  (4.57), one establishes that (4.53b) is an identity.

We are now ready to prove Proposition 4.3 in the case where $i=N-1$.
\begin{proposition}
In $N$ dimensions
\begin{equation}
\begin{split}
&\Res [H_N,B_{N-1}^+]=\\
&\qquad \qquad \Res \biggl\{ B_{N-1}^+ \Bigl( 2\sum_{i=1}^N D_i +(N-1)(1+\beta)
\Bigr) +G_N \tilde N_N  \biggr\}
\end{split}
\end{equation}
where $G_N$ are some unspecified expressions.
\end{proposition}
Proof.  Recall that $H_N= \sum_{i=1}^N D_i^2$.  Since each term of $B_{N-1}^+$
has a product of $N-1$ distinct variables as factor on the left,
it follows from Lemma  4.10 that the commutator of $H_N$ and $B_{N-1}^+$ is of
the form
\begin{equation}
[H_N,B_{N-1}^+]= \sum_{j=1}^{N} z_1 \dots  \hat z_j  \dots  z_{N} \mathcal F_j
\end{equation}
with $\mathcal F_j$ expressions involving the $D_{\ell}$ and $K_{mn}$.  Since
$H_N$ and $B_{N-1}^+$ are invariant under the action of the symmetric group
$S_N$,
 an argument similar to the one given in the proof of Proposition 4.8 shows
that
to compute the commutator $[H_N,B_{N-1}^+]$, one only needs to determine the
terms with $z_2 \dots z_N$ as factor on the left and to symmetrize the result.
This
is the approach that we shall take.

We initiate the induction by proving (4.61) when $N=2$.  In this special case
 we can write
\begin{equation}
H = \Res (D_1^2+ D_2^2) = \Res \Bigl\{ (D_1+D_2)(D_1+D_2+\beta)-2\tilde N_2
\Bigr\}.
\end{equation}
It is then easy to see with the help say, of Proposition 4.8 and using
$\Res \sum_i D_i $  =  $\sum_i z_i \frac{\partial}{\partial z_i}$ that
\begin{equation}
\Res [H_2,B_1^+]=\Res \Bigl\{ B_1^+ \bigl( 2(D_1+D_2) + 1+\beta) \bigr)
+G_2 \tilde N_2  \Bigr\}.
\end{equation}
We shall now show that (4.61) is satisfied in $N$ dimensions if it is obeyed
in dimensions lower than $N$.  To this end, we shall need the following
expression for $B_{N-1}^+$:
\begin{equation}
B_{N-1}^+ = z_N B^+_{N-2,\{1,\dots,N-1\}} \bigl( D_N +(N-1)\beta \bigr) +
z_1 \dots z_{N-1} D_{1,\{1,\dots,N-1\}}
\end{equation}
in which the dependance on $z_N$ is isolated.  In the notation (4.52),  we
shall also write
\begin{equation}
H_N = H_{N-1} + D_N^2,\qquad H= \Res  H_N.
\end{equation}
With the help of Lemma  4.10, this leads to
\begin{equation}
\begin{split}
& \Res [H_N,B_{N-1}^+] \bigm|_{z_1 \sim 0}=\\
& \qquad \qquad \Res \biggl\{ [H_{N-1}, z_N B^+_{N-2,\{1,\dots,N-1\}}] \bigl(
D_N+(N-1)\beta \bigr) \\
&\qquad \qquad ~~~~+z_2 \dots z_N D_{1,\{2,\dots,N-1\}}[H_{N-1},D_N] +
[D_N^2,z_2 \dots z_N D_{1,\{2,\dots,N\}}] \\
&\qquad \qquad \qquad \qquad \qquad \qquad \qquad \quad+[D_N^2,z_1 \dots
z_{N-1}] D_{1,\{1,\dots,N-1\}} \biggr\}
\biggm|_{z_1 \sim 0 }
\end{split}
\end{equation}
We now invoke  Corollary 4.11 and the definition of $B_{N-2}^+$ to assert
that $[H_{N-1},z_N B^+_{N-2,\{1,\dots,N-1\}}]\bigm|_{z_1 \sim 0}$ and
$z_N [H_{N-1}^{\{N-1\}}, B^{+\{N-1\}}_{N-2}]\bigm|_{z_1 \sim 0}$ have the same
symbol structure.  This observation and the fact that $\bigl( D_N +(N-1) \beta
\bigr)$
is invariant under $S_{N-1}$ allow one to compute the first term in the
right-hand side of (4.67) from the induction hypothesis.  Indeed we have
\begin{equation}
\begin{split}
&\Res^{\{N-1\}} [H_{N-1},z_N B^+_{N-2,\{1,\dots,N-1\}}]  \bigm|_{z_1 \sim 0}\\
&\qquad  = \Res^{\{N-1\}} \Biggl[ z_N \biggl\{ z_2 \dots z_{N-1}
D_{1,\{2,\dots,N-1\}} \\
& \qquad \qquad  \times \Bigl(2\sum_{i=1}^{N-1} D_i + (N-2)+(N-2)\beta \Bigr) +
G_{N-1} \tilde N_{N-1}  \biggr\}
\Biggr] \Biggm|_{z_1 \sim 0}
\end{split}
\end{equation}
where it is understood that the Dunkl operators depend on the variables
$z_1,\dots,$  $z_N$.

The other terms on the right-hand side of (4.67) are computed with the help of
Lemma  4.12 (with $\ell = 2$) and Lemma  4.10.  Putting everything together and
using the fact that $\tilde N_{N-1} \bigl( D_N+(N-1)\beta \bigr)=\tilde N_N$,
one
sees that (4.67) reduces to
\begin{equation}
\begin{split}
& \Res [H_N,B_{N-1}^+] \bigm|_{z_1 \sim 0}=\\
&~ z_2 \dots z_N \Res \biggl\{ D_{1,\{2,\dots,N-1\}} \Bigl(2\sum_{i=1}^{N-1}
D_i + (N-2)(1+\beta) \Bigr) \bigl( D_N+(N-1)\beta \bigr)\\
&~~~~~~\qquad+ \Bigl( 2D_N + 1+\beta -
\beta(1-K_{1N})(D_N-\beta)-2\beta(1-K_{1N})
\Bigr) D_{1,\{2,\dots,N\}}\\
&~~~~~~~~~~~~\qquad \qquad \qquad \qquad \qquad+ \beta D_{1,\{2,\dots,N-1\}}
(D_N^2-D_1^2) + G^{\prime}_N
\tilde N_N \biggr\} \biggm|_{z_1 \sim 0}
\end{split}
\end{equation}
where $G^{\prime}_N$ is another unspecified quantity.

The first two terms between the curly brackets on the right-hand side of (4.69)
can be rewritten in the form
\begin{equation}
\begin{split}
&   \biggl\{ D_{1,\{2,\dots,N\}} \Bigl( 2\sum_{i=1}^{N-1} D_i
+(N-2)(1+\beta) \Bigr)+D_{1,\{2,\dots,N-1\}}[2\sum_{i=1}^{N-1} D_i, D_N] \\
& \qquad \qquad  +D_{1,\{2,\dots,N-1\}} \bigl( 2 D_N + 1 + \beta \bigr)
+ 2 [D_N, D_{1,\{2,\dots,N\}}] \\
& \qquad \qquad \qquad  + \Bigl( - \beta(1-K_{1N})(D_N-\beta)-2\beta(1-K_{1N})
\Bigr) D_{1,\{2,\dots,N\}} \biggr\}.
\end{split}
\end{equation}
Using formulas (4.53a) and (4.53b) with $ \ell = 1$ in (4.70), one can recast
(4.69) as follows:
\begin{equation}
\begin{split}
& \Res [H_N,B_{N-1}^+] \bigm|_{z_1 \sim 0}= \\
& \Res z_2 \dots z_N \biggl\{ D_{1,\{2,\dots,N\}} \Bigl(2\sum_{i=1}^{N} D_i +
(N-1)(1+\beta) \Bigr)\\
&-\beta(1-K_{1N})(D_N - \beta) D_{1,\{2,\dots,N\}} + \beta
D_{1,\{2,\dots,N-1\}} (D_N^2-D_1^2) + G^{\prime}_N
\tilde N_N \biggr\} \biggm|_{z_1 \sim 0}.
\end{split}
\end{equation}
With the help of (4.53b) with $\ell = 1$ again, we find
\begin{equation}
\begin{split}
& \Res \Bigl( -\beta (1-K_{1N})D_N D_{1,\{2,\dots,N\}} \Bigr) =\\
&~~~~~~~~~~~~~~~\Res \Bigl( -\beta D_{1,\{2,\dots,N-1\}} (D_N^2-D_1^2)
-\beta^2 (1-K_{1N}) D_{1,\{2,\dots,N\}} \Bigr)
\end{split}
\end{equation}
so that (4.71) becomes
\begin{equation}
\begin{split}
& \Res [H_N,B_{N-1}^+] \bigm|_{z_1 \sim 0}= \\
& \quad  \Res \biggl\{
 z_2 \dots z_N  D_{1,\{2,\dots,N\}} \Bigl( 2\sum_{i=1}^{N} D_i^2
+(N-1)(1+\beta) \Bigr) +G_N'' \tilde N_N  \biggr\} \biggm|_{z_1 \sim 0}.
\end{split}
\end{equation}
Upon symmetrizing the right-hand side of (4.73), we finally obtain
\begin{equation}
\Res [H_N,B_{N-1}^+] = \Res \biggl\{
 B_{N-1}^+  \Bigl( 2\sum_{i=1}^{N} D_i^2
+(N-1)(1+\beta) \Bigr) + G_N \tilde N_N \biggr\},
\end{equation}
thereby proving Proposition 4.13.

Given that Proposition 4.3 is true when $i=N-1$, we shall fix $i$ to be
equal to $N-1$ and  extend the number of variables from $N$ to an arbitrary
larger number $N'$.  With
the symmetry of $H_{N'}$ and $B_{N-1}^+$ under the exchange of
the variables $ z_1,\dots,z_{N'}$  allowing to insert
$\Res$ (meaning $\Res^{\{N'\}}$) in front of every operator, Proposition~4.3
will follow from the next proposition.
\begin{proposition}
In $N'$ dimensions,
\begin{equation}
\begin{split}
& [H, \Res B_{N-1}^+]  = \Res B_{N-1}^+ \Bigl\{ 2 \Res \sum_{i=1}^{N'} D_i
+(N-1) + \beta (N-1)(N'-N+1) \Bigr\}\\
& \qquad \qquad \qquad \qquad \qquad \qquad \qquad  \qquad + \Res
\sum_{\begin{subarray}{c} J \subset \{1,\dots,N'\}\\ |J|=N-1 \end{subarray}}
G_{N',J} \tilde N_{N,J},
\end{split}
\end{equation}
with $G_{N',J}$ an unspecified expression.
\end{proposition}
In order to prove (4.75), we shall need the next three lemmas
\begin{lemma}
In $N'$ dimensions, one has
\begin{equation}
\begin{split}
& \Res D_{k,J} =  \Bigl( z_{j_1} \frac{\partial}{\partial z_{j_1}}
+ \beta \sum_{i \in \{2,\dots,{\ell}\}} \frac{z_{j_1}}{z_{j_1} - z_{j_i}}
(1-K_{j_1 j_i}) + k\beta \Bigr) \dots \\
& \qquad \qquad \qquad \qquad \qquad \dots \Bigl( z_{j_{\ell}}
\frac{\partial}{\partial z_{j_{\ell}}} +(k+\ell-1)\beta \Bigr)=\Res^J
D_{k,J}^{(J)}\\
& \qquad \qquad \qquad \qquad \qquad k=0,1,2,\dots
\end{split}
\end{equation}
for all subsets $J \subset \{1,\dots,N'\}$ of cardinality $|J|=\ell$
such that, $j_{\kappa}<j_{\kappa+1}$ if $j_{\kappa} \in J$ and $1 \leq \kappa <
\ell$.
\end{lemma}
Proof. This lemma follows from the definition of $D_{k,J}$ and the fact
that $K_{ij}D_{\ell}=D_{\ell}K_{ij}$ if $\ell \neq i,j$.  It expresses the fact
that $\Res D_{k,J}$ only depends upon the variables $z_i$, $i \in J$.

\begin{corollary} For all sets $J$ defined in Lemma~4.15, $\Res B_{N-1,J}^+$
only depends upon the variables $z_i$, $i \in J$.
\end{corollary}
Let  $p_i = z_i \partial/\partial z_i$ and
$
A_{ij}
	=  \frac{z_i+z_j}{z_i - z_j} \Bigl( z_i \frac{\partial}{\partial z_i} - z_j
\frac{\partial}{\partial z_j} \Bigr)$, and recall that the Hamiltonian reads
\begin{equation}
H^{(J)}
	= \sum_{i \in J} p_i^2 + \beta \sum_{\begin{subarray}{c} i<j\\i,j \in J
\end{subarray}} A_{ij},
\end{equation}
 in terms of the variables $z_i, i \in J$.  This expression appears in the
following decomposition of the commutator $[H, \Res B_{N-1}^+]$.

\begin{lemma} For any $N' > N$,
\begin{equation}
[H, \Res B_{N-1}^+]
	= \sum_{m=1}^{N'+1 - N} (-)^{m+1} \sum_{\begin{subarray}{l} J \subset
\{ 1, \dots , N'\}\\ |J| = N' - m \end{subarray}} [H^{(J)}, \Res B_{N-1,J}^+],
\end{equation}
with $B_{N-1,J}^+ \equiv z_J D_{1,J}$  if  $|J| = N-1$.
\end{lemma}

Proof.  The proof  is combinatorial.  We need to show that the right-hand side
of (4.78) contains each summand of
$\Res B_{N-1}^+$ commuted  once and only once  with
all the  parts of $H = H^{\{1,\dots , N'\}}$ given in (4.77).  In view of the
$S_{N'}$-symmetry,  it is
sufficient to show  that this is true for one summand say,  $\Res z_1 \dots
z_{N-1}
D_{1,\{ 1,\dots , N-1\}}$.  From Lemma~4.15, this operator depends only on the
variables  $z_1, \dots , z_{N-1}$ and it thus suffices to look for the number
of times $\Res z_1 \dots z_{N-1} D_{1,\{1, \dots, N-1\}}$  is commuted with
$p_i^2$, $A_{ij}$ and $A_{ik}$ ($i,j \in \{ 1, \dots , N-1\}$,
$~ k \notin \{ 1,\dots , N-1\}$).   Since the sets $J$ of (4.78) must in this
case contain $\{1,\dots,N-1\}$,
we  see that the terms of the right-hand side of this equation that involve
$\Res z_1 \dots
z_{N-1} D_{1,\{1,\dots , N-1\}}$ are
\begin{equation}
\sum^{N'+1-N}_{m=1} (-)^{m+1} \sum_{\begin{subarray}{c} J' \subset
\{ N,\dots , N'\}\\|J'| =  N'+1-N-m \end{subarray}} [H^{\{1,\dots, N-1\}\cup
J'},
z_1 \dots  z_{N-1} \Res D_{1,\{1,\dots , N-1\}}].
\end{equation}
Each $H^{\{1,\dots N-1\}\cup J'}$ contains exactly one $p_i^2$ and one
$A_{ij}$.   For $|J'| = N' + 1 - N-m$, both $p_i^2$ and $A_{ij}$ are seen to
appear $\binom{N'+1-N}{N'+1-N-m}$= $\binom{N'+1-N}{m}$ times in the commutators
of (4.79), thus for a total of $\sum_{m=1}^{N'+1-N}(-)^{m+1}
\binom{N'+1-N}{m}=1$ occurence, as required.  In the case of the $A_{ik}$ ,
since $k$ must be present in $J'$, it  shows up  $\binom{N'-N}{N'-N-m}=
\binom{N'-N}{m}$ times in (4.79),  if $m<N'+1-N$ and does not appear otherwise.
 Performing the sum over $m$, we also find  $\sum_{m=1}^{N'-N}(-)^{m+1}
\binom{N'-N}{m}=1$, which proves the lemma.  We see from the last summation
that (4.78) does not hold  for $N'=N$, since, in this case,  $A_{iN}$  never
appears.

Let us finally prove.
\begin{lemma}
For all sets $J$ defined in Lemma  4.15,
\begin{equation}
[H^{(J)}, \Res  z_J D_{1,J}]= z_J D_{1,J} \Bigl( 2 \Res
\sum_{i \in J} D_i + \ell \Bigr),
\end{equation}
with $\ell=|J|$.
\end{lemma}
Proof.  We have from (4.4) and Lemma  4.15 that
\begin{equation}
[H^{(J)} ,  \Res  D_{1,J}]=
 \Res^J [H^{(J)}_{\ell} ,   D_{1,J}^{(J)}]=0,
\end{equation}
 and from Lemma  4.10 that
\begin{equation}
[H^{(J)},z_J] = \Res^J [H^{(J)}_{\ell},z_J] = z_J (2 \Res^J
\sum_{i \in J} D_i^{(J)}+ \ell) = z_J (2 \Res \sum_{i \in J} D_i + \ell).
\end{equation}
Using (4.4) and Lemma  4.15 again to commute $\Res \sum_{i \in J} D_i$ and
$\Res D_{1,J}$,
(4.80) is seen to hold.

  We are now ready to prove Proposition 4.14.  From Lemma~4.18 and
Proposition~4.3, Proposition~4.14 is seen to be true for the cases  $N'=N-1$
and $N$.  Remarkably, when $N'> N$, Lemma  4.17
shows that $[H, \Res B_{N-1}^+]$ can be decomposed in commutators involving
less than $N'$ variables.  The induction process is thus greatly  simplified.
We use Corollary~4.16 to identify  $\Res^{J} B_{N-1,J}^{+(J)}$ and  $\Res
B_{N-1,J}^+$ and take, by hypothesis,
\begin{equation}
\begin{split}
& [H^{(J)}, \Res  B_{N-1,J}^+]= \Res B_{N-1,J}^+ \Bigl\{ 2 \Res
\sum_{i \in J} D_i + (N-1) \\
& \qquad \qquad + \beta (N-1)(N'+1-m-N) \Bigr\} + \Res G_{N'-m,J} \tilde
N_{N-1,J},
\end{split}
\end{equation}
for $m=1,\dots,N'+1-N$, with $|J|=N'-m$, in (4.78).   There thus remain to be
shown that (4.75) is obtained upon  performing the sum in (4.78).

Inserting (4.83) in (4.78) gives
\begin{equation}
\begin{split}
&[H, \Res B_{N-1}^+]=\sum_{m=1}^{N'+1-N} (-)^{m+1}\sum_{\begin{subarray}{c}
J \subset \{1,\dots,N'\}\\|J|=N'-m \end{subarray}} \biggl\{ \Res B_{N-1,J}^+
\Bigl( 2 \Res  \sum_{i \in J} D_i \\
&\qquad + (N-1)
 + \beta (N-1)(N'+1-m-N) \Bigr) + \Res G_{N'-m,J} \tilde N_{N-1,J} \biggr\}.
\end{split}
\end{equation}
In order to compute the right-hand side, as in the proof of Lemma  4.17, we
look at the part
 of $\Res B_{N-1}^+$ involving $z_1 \dots z_{N-1} \Res
 D_{1,\{1,\dots,N-1\}}$ in (4.84).
Considering the terms that  appear on the left of this particular term in
(4.84), we get
\begin{equation}
\begin{split}
& \sum_{m=1}^{N'+1-N}  (-)^{m+1}
 \sum_{\begin{subarray}{c} J' \subset \{N,
\dots,N'\}\\|J'|=N'+1-N-m \end{subarray}}  \Bigl( 2 \Res \sum_{i \in
\{1,\dots,N-1\}
\cup J'}
 D_i \\
& \qquad \qquad \qquad \qquad\qquad \qquad+ (N-1)
 + \beta (N-1)(N'+1-m-N) \Bigr) \\
&  = \sum_{m=1}^{N'+1-N} (-)^{m+1} \binom{N'+1-N}{m}  \Bigl( 2
\Res \sum_{i =1}^{N-1} D_i
+ (N-1)\\
& \qquad \quad + \beta (N-1)(N'+1-m-N) \Bigr)
+ \sum_{m=1}^{N'-N} (-)^{m+1}
\binom{N'-N}{m}   2 \Res  \sum_{i =N}^{N'} D_i \\
&   = 2 \Res  \sum_{i=1}^{N'} D_i +(N-1)+\beta (N-1)(N'+1-N),
\end{split}
\end{equation}
where we used
\begin{equation}
 \sum_{m=1}^{N'+1-N} (-)^{m+1} \binom{N'+1-N}{m} m \propto \sum_{m=0}^{N'-N}
(-)^{m} \binom{N'-N}{m} = 0.
\end{equation}
Since the result is the same for every part of $ B_{N-1}^+$, setting
$\sum_{m=1}^{N'+1-N} (-)^{m+1} \* \Res G_{N'-m,J}$ = $\Res  G_{N',J}$, we get
the desired result, namely that
\begin{equation}
\begin{split}
& [H,\Res B_{N-1}^+] = \Res B_{N-1}^+ \Bigl\{ 2 \Res \sum_{i=1}^{N'} D_i
+(N-1) + \beta (N-1)(N'-N+1) \Bigr\}\\
& \qquad \qquad \qquad \qquad \qquad \qquad \qquad + \Res
\sum_{\begin{subarray}{c} J \subset \{1,\dots,N'\}\\ |J|=N-1 \end{subarray}}
G_{N',J} \tilde N_{N,J},
\end{split}
\end{equation}
thereby proving  Proposition 4.14 and, as a consequence, Proposition 4.3.

\subsection{The coefficient $c_{\lambda}$}
This last subsection is devoted to the proof of the expression (3.14)-(3.15) of
the coefficient $c_{\lambda}$ of Theorem~3.1.
\begin{proposition}
For any partition $\lambda=(\lambda_1,\dots,\lambda_{\ell})$,
\begin{equation}
B_{\ell}^+ m_{\lambda}= a_{\lambda} m_{\lambda+1} + \sum_{\lambda^{\prime} <
\lambda+1} d_{\lambda,\lambda'} m_{\lambda^{\prime}}
\end{equation}
with $a_{\lambda}=(\lambda_{\ell}+\beta)\dots (\lambda_1+\ell \beta)$,
$\lambda+1=(\lambda_1+1,\dots,\lambda_{\ell}+1)$ and the
  $d_{\lambda,\lambda'}$'s certain coefficients.
\end{proposition}
By the triangularity of the Jack polynomials on the monomial basis
$\{m_{\rho}\}$,
this proposition implies that $ B_{\ell}^+ J_{\lambda}=a_{\lambda}
J_{\lambda+1}=(\lambda_{\ell}+\beta)\dots (\lambda_1+\ell \beta)J_{\lambda+1}$,
which, by successive applications of the creation operators, gives the form
(3.14) of $c_{\lambda}$.

The proof of  Proposition 4.19 is done by looking for the coefficient appearing
in (4.88) in front of the  term $z_1^{\lambda_{\ell}+1} \dots
z_{\ell}^{\lambda_1+1}$ of $m_{\lambda+1}$. We shall see that this term can  be
generated through simple operations only.

 Since $z_1^{\lambda_{\ell}+1} \dots z_{\ell}^{\lambda_1+1}$  involves only the
variables $z_1,\dots,z_{\ell}$, it can only be generated in (4.88) from the
part $z_1 \dots z_{\ell} D_{1,\{1,\dots,\ell\}}$ of $B_{\ell}^+$.  Moreover,
\begin{equation}
m_{\lambda}= \sum_{{\text{distinct}}} z_1^{\lambda_{p(1)}} \dots
z_{\ell}^{\lambda_{p(\ell)}} + {\text {terms involving at least one other
variable}}
\end{equation}
and  $\Res D_{1,\{1,\dots,\ell\}}$, from Lemma~4.15 , depends only on the
variables
$z_1,\dots,z_{\ell}$.  Thus, only the action of $z_1 \dots z_{\ell}
\Res D_{1,\{1,\dots,\ell\}}$ on $  \sum_{{\text{distinct}}}
z_1^{\lambda_{p(1)}} \dots z_{\ell}^{\lambda_{p(\ell)}}$ can generate
$z_1^{\lambda_{\ell}+1} \dots z_{\ell}^{\lambda_1+1}$.

We now establish the following lemma.
\begin{lemma}
For any $i<j$ and $P \in S_{\ell}$, the expansion of
\begin{equation}
\frac{z_i}{z_i-z_j} (1-K_{ij})~~~  z_1^{\lambda_{p(1)}} \dots
z_{\ell}^{\lambda_{p(\ell)}}
\end{equation}
does not contain any term in $z_1^{\lambda_{\ell}} \dots z_{\ell}^{\lambda_1}$,
 where $\lambda_1 \geq \lambda_2 \geq \dots \geq \lambda_{\ell}$.
\end{lemma}
Proof.  $\forall i<j$
\begin{equation}
\begin{split}
&\frac{z_i}{z_i-z_j} (1-K_{ij})~~~  z_1^{\lambda_{p(1)}}
\dots z_{\ell}^{\lambda_{p(\ell)}} =  z_1^{\lambda_{p(1)}} \dots \hat z_{i}
\dots \hat z_j \dots z_{\ell}^{\lambda_{p(\ell)}} \\
& \quad \times {\cases  -z_i^{\lambda_{p(j)}} z_j^{\lambda_{p(i)}} -
\sum\limits_{k=0}^{\lambda_{p(j)}-\lambda_{p(i)}-2} z_i^{\lambda_{p(j)}-1-k}
z_j^{\lambda_{p(i)}+1+k},&{\text{if}}~ \lambda_{p(j)} > \lambda_{p(i)}\\
 z_i^{\lambda_{p(i)}} z_j^{\lambda_{p(j)}} +
\sum\limits_{k=0}^{\lambda_{p(i)}-\lambda_{p(j)}-2} z_i^{\lambda_{p(i)}-1-k}
z_j^{\lambda_{p(j)}+1+k} ,&{\text{if}}~ \lambda_{p(j)} < \lambda_{p(i)}\\
 0,   &{\text{if}}~ \lambda_{p(j)} = \lambda_{p(i)}. \endcases}
\end{split}
\end{equation}
Hence, the term associated to the partition $(\lambda_1,\dots,\lambda_{\ell})$
in the right-hand side of (4.91)
is
\begin{equation}
\cases  -  z_1^{\lambda_{p(1)}} \dots z_{i}^{\lambda_{p(j)}}
\dots z_j^{\lambda_{p(i)}} \dots z_{\ell}^{\lambda_{p(\ell)}} \qquad
&{\text{if}}~ \lambda_{p(j)} > \lambda_{p(i)}\\
  z_1^{\lambda_{p(1)}} \dots z_{i}^{\lambda_{p(i)}} \dots z_j^{\lambda_{p(j)}}
\dots z_{\ell}^{\lambda_{p(\ell)}} &{\text{if}}~ \lambda_{p(j)} <
\lambda_{p(i)}\\
 0 &{\text{if}}~ \lambda_{p(j)} = \lambda_{p(i)},
\endcases
\end{equation}
from where we see that the exponent of $z_j$ is always less than the one of
$z_i$ and  that, consequently, the term $z_1^{\lambda_{\ell}} \dots
z_{\ell}^{\lambda_1}$ cannot appear.

{}From Lemma  4.20, Lemma  4.15 in the case $J=\{1,\dots,\ell\}$, and the
triangularity of the action of $\frac{z_i}{z_i-z_j} (1-K_{ij})$ on $
z_1^{\lambda_{p(1)}} \dots z_{\ell}^{\lambda_{p(\ell)}}$, we finally find that
$z_1^{\lambda_{\ell}+1} \dots z_{\ell}^{\lambda_1+1}$ can only appear in the
following way through the action of $z_1 \dots z_{\ell}
\Res^{\{N\}}D_{1,\{1,\dots,\ell\}}$ on $m_{\lambda}$
\begin{equation}
\begin{split}
& z_1 \dots z_{\ell} (z_1 \frac{\partial}{\partial z_1} + \beta)\dots (z_{\ell}
\frac{\partial}{\partial z_{\ell}} + \ell \beta)~~ z_1^{\lambda_{\ell}} \dots
z_{\ell}^{\lambda_1} \\
& \qquad \qquad \qquad \qquad \qquad \qquad \qquad =
(\lambda_{\ell}+\beta)\dots(\lambda_1+\ell \beta)~~ z_1^{\lambda_{\ell}+1}
\dots z_{\ell}^{\lambda_1+1}\\
& \qquad \qquad \qquad \qquad \qquad \qquad \qquad= a_{\lambda}~
z_1^{\lambda_{\ell}+1} \dots z_{\ell}^{\lambda_1+1}.
\end{split}
\end{equation}
This proves Proposition 4.19.

We may remark that the coefficient $c_{\lambda}$ appears in (3.13) only because
we use
Jack polynomials that are monic.  Stanley uses instead in \cite{8} the
normalization which is defined by taking $v_{\lambda,(1^n)}=n !$ if
$|\lambda|=n$ in the expansion (2.18) of the Jack polynomials $J_{\lambda}$ in
terms of the symmetric monomials.  We note that in this normalization
$c_{\lambda}=\beta^n  v_{\lambda,\lambda}$ therefore, had we used the
normalization of Stanley and redefined the creation operators according to
$B_i^+ \to 1/\beta^i B_i^+$, we would have found \cite{18} that
$$ J_{\lambda}(z;1/\beta)
	= \bigl( B_{N-1}^{+}\bigr)^{\lambda_{N-1}} \dots
\bigl( B_{2}^{+}\bigr)^{\lambda_2-\lambda_3}\bigl( B_{1}^{+}\bigr)^
{\lambda_{1}-\lambda_2} \cdot 1 $$
without any proportionality constant.
\section{Conclusion}
Our general objective is to develop a completely algebraic treatment of the
Calogero-Sutherland model.  We would thus hope to identify the abstract
structure of its full dynamical algebra and to work out its relevant
representations.  This should in principle allow one to obtain algebraically
all physically interesting quantities.

We believe that the results presented in this paper provide important clues
toward the resolution of these questions.  In fact, they have already allowed
us to make progress and to formulate in this connection remarkable conjectures.
 Clearly, the creation operators $B_i^+$, $i=1,\dots,N$, should be among the
generators of the dynamical algebra.  One would thus want to know their action
on the wave functions of the CS model, that is on the Jack polynomials.
However, since these operators do not commute among themselves and enter in a
definite
order in formula (3.13)(or (3.17)) for the excited wave functions, it is not
straightforward to obtain these actions.  We have nevertheless a conjecture
for this.  We also have expressions for annihilation operators $B_i^-$,
$i=1,\dots,N$ and a similar conjecture giving their action on the wave
functions.
  These developments allow us to evaluate in particular the norm of the
Jack polynomials in an algebraic fashion.  We shall report on these results
 in a forthcoming publication \cite{17}.  A major problem that remains however
is to determine the structure relations that the creation and annihilation
operators obey.

The results of this paper also have mathematical applications.  Various
conjectures involving the Jack polynomials have been made.  It turns out
\cite{18} that
 formula (3.13) readily implies a weak form of a famous conjecture due to
 Macdonald and Stanley, namely that in the normalization of Stanley the
$v_{\lambda \mu}$ in (2.18b) are polynomials in $\beta^{-1}$ with integer
coefficients (see the remark at the end of section~4.3).  We believe that this
Rodrigues formula that we have obtained will provide a useful tool to further
advance the proofs of the outstanding conjectures on the Jack polynomials.

Finally, it is expected that the results presented here extend to the
relativistic
 generalization of the Calogero-Sutherland model \cite{19}.  This requires a
$q$-deformation of  our constructs and should yield a Rodrigues-type formula
for the Macdonald polynomials.  It is also of interest to consider in the same
vein models and special functions associated to root lattices other than
$A_{N-1}$.  We are currently studying these questions and hope to report on
them in the near future.

\begin{acknow}
We are grateful to Fran\c{c}ois Bergeron, Adriano Garsia and Jean Le  Tourneux
for comments and discussions.

\noindent This work has been supported in part through funds provided by NSERC
(Canada) and FCAR (Qu\'ebec).  One of us (L.L.) holds a NSERC postgraduate
scholarship.
\end{acknow}

\end{document}